\documentclass[submission]{eptcs}
\providecommand{\event}{ACT 2020}

\newcommand{\PartFn}{\mathsf{PartFn}}
\newcommand{\Poset}{\mathsf{Poset}}
\newcommand{\Set}{\mathsf{Set}}
\newcommand{\Rel}{\mathsf{Rel}}
\newcommand{\Geom}{\mathsf{Geom}}
\newcommand{\Data}{\mathsf{Data}}
\newcommand{\Sem}{\mathsf{Sem}}
\newcommand{\Bool}{\mathsf{Bool}}
\newcommand{\BoolAlg}{\mathsf{BoolAlg}}

\newcommand{\MoveTo}{\textbf{MoveTo}}
\newcommand{\Pick}{\textbf{Pick}}
\newcommand{\Place}{\textbf{Place}}

\renewcommand{\AA}{\mathcal{A}}
\newcommand{\BB}{\mathcal{B}}
\newcommand{\CC}{\mathcal{C}}
\newcommand{\DD}{\mathcal{D}}

\newcommand{\from}{\leftarrow}

\newcommand{\epi}{\twoheadrightarrow}

\newcommand{\id}{\mathrm{id}}
\newcommand{\el}{\mathrm{el}}

\newcommand{\R}{\mathbb{R}}

\newcommand{\SE}{\mathsf{SE}}

\newcommand{\Atom}{\mathsf{Atom}}
\newcommand{\List}{\mathsf{List}}
\newcommand{\Pt}{\mathsf{Pt}}

\newcommand{\arity}{\mathsf{ar}}
\newcommand{\tp}{\mathsf{tp}}

\newcommand{\pbcorner}[1][rd]{\ar@{}[rd]|(.3){\scalebox{2}{$\lrcorner$}}}

\RequirePackage{array}

\usepackage{underscore}

\usepackage{amsmath,amssymb,amsthm}
\usepackage[all,cmtip,2cell]{xy}
	\UseTwocells

\usepackage{ifthen}
\usepackage{mathtools}
\usepackage{hyperref}

\usepackage{tikz}
\usetikzlibrary{positioning, calc, decorations.markings}

\usepackage{tikz-cd}

\theoremstyle{plain}
\newtheorem{prop}{Proposition}

\theoremstyle{definition}
\newtheorem{defn}{Definition}
\newtheorem{ex}{Example}

\DeclareMathOperator{\Ob}{Ob}

\title{Symmetric Monoidal Categories with Attributes}

\author{Spencer Breiner
\institute{National Institute of Standards and Technology \\
Gaithersburg, MD, USA}
\email{Spencer.Breiner@nist.gov}
\and
John S.\ Nolan\thanks{Corresponding author.}
\institute{University of Maryland \\
College Park, MD, USA}
\email{jnolan13@terpmail.umd.edu}
}
\def\titlerunning{Symmetric Monoidal Categories with Attributes}
\def\authorrunning{S.\ Breiner and J.\ S.\ Nolan}

\begin{document}

\maketitle

\begin{abstract}
When designing plans in engineering, it is often necessary to consider attributes associated to objects, e.g.\ the location of a robot.
Our aim in this paper is to incorporate attributes into existing categorical formalisms for planning, namely those based on symmetric monoidal categories and string diagrams.
To accomplish this, we define a notion of a ``symmetric monoidal category with attributes.''
This is a symmetric monoidal category in which objects are equipped with retrievable information and where the interactions between objects and information are governed by an ``attribute structure.''
We discuss examples and semantics of such categories in the context of robotics to illustrate our definition.
\end{abstract}

\section{Introduction}

% preliminaries: assume familiarity with SMCs

Symmetric monoidal categories (SMCs) and the related graphical syntax of string diagrams have recently been used to great effect in representing and understanding plans and processes (see e.g.\ \cite{CoeckeKissinger, sevensketches}).
In large part, this is because the syntax of SMCs enables users to reasonably interpret objects in a SMC as \emph{resources} and morphisms in a SMC as \emph{processes}.

Here we extend this interpretation to include a distinction between physical and informational resources, \emph{entities} and \emph{values} or \emph{data}, which behave quite differently. Our approach is inspired by the relationship between classical and quantum information in categorical quantum mechanics \cite{CoeckeKissinger}. We are motivated by foundational concerns in a related paper \cite{robotspaper}, where we use this framework to link high-level action planning and low-level path planning in the context of robotics.
	
As an example, consider the following equation (where the diagrams are to be read top-to-bottom), which axiomatizes one of the operations of a robot arm:
% =======================================================================================================================
% SHRINK ME WHEN NAVIGATING THE PAPER
% =======================================================================================================================
	\begin{equation}\begin{tabular}{ccc}
			\begin{tabular}{c}
				\begin{tikzpicture}
				\tikzset{
	% Nodes
	AttrR/.style = {font={\Large$\rhd$}},
	AttrL/.style = {font={\Large$\lhd$}},
	MoveTo/.style = {draw,font={\MoveTo}},
	Pick/.style = {draw,font={\Pick}},
	Place/.style = {draw,font={\Place}},
	Circ/.style = {draw,circle,scale=0.75},
	Rect/.style = {draw,rectangle},
	%Edges
	ob/.style = {out=-90, in=90, looseness=1},
	info/.style = {dashed,out=-90, in=90, looseness=1},
	attrR/.style = {dashed,out=0, in=90, looseness=1},
	attrL/.style = {dashed,out=180, in=90, looseness=1},
	attrRL/.style = {dashed,out=0, in=180, looseness=1},
	attrLR/.style = {dashed,out=180, in=0, looseness=1},
%	attr/.style = {dashed,out=0, in=90, looseness=1},
%	urattr/.style = {dashed,out=-90, in=180, looseness=1},
%	ldattr/.style = {dashed,out=180, in=90, looseness=1},
%	ulattr/.style = {dashed,out=-90, in=0, looseness=1},
%	lrattr/.style = {dashed,out=180, in=0, looseness=1,shorten >=-8pt, shorten <=-8pt},
}

\pgfmathsetmacro{\port}{3mm}

				% Frame
				\coordinate (NW) at (0,2); \coordinate (NE) at (4,2); 
				\coordinate (SW) at (0,0); \coordinate (SE) at (4,0);
				%\draw[black] (NW) rectangle (SE);
				
				%Nodes
				\node[MoveTo] (MoveTo1) at ({$(NW)!0.5!(NE)$} |- {$(SW)!0.5!(NW)$}) {};
				\node[Circ] (Circ1) at ({$(NW)!0.66!(NE)$} |- {$(SW)!0.9!(NW)$}) {};
				\node[AttrR] (AttrR1) at ({$(NW)!0.54!(NE)$} |- {$(SW)!0.1!(NW)$}) {};
				
				% Autocreate port names
				\pgfmathsetmacro{\diagramIn}{2}
				\pgfmathsetmacro{\diagramOut}{1}
				\pgfmathsetmacro{\nAttrL}{0}
				\pgfmathsetmacro{\nAttrR}{1}
				\pgfmathsetmacro{\nMoveTo}{1}
				\pgfmathsetmacro{\nPick}{0}
				\pgfmathsetmacro{\nPlace}{0}
				\foreach \x [evaluate=\x as \dx using \x/(\diagramIn + 1)] in {1,...,\diagramIn}{
	\coordinate (in\x) at ([yshift=\port]$(NW)!\dx!(NE)$);
}

\foreach \x [evaluate=\x as \dx using \x/(\diagramOut + 1)] in {1,...,\diagramOut}{
	\coordinate (out\x) at ([yshift=-\port]$(SW)!\dx!(SE)$);
}

\ifthenelse	{\nAttrR > 0}
			{\foreach \i in {1,...,\nAttrR}{
				\coordinate (AttrR\i in) at ([xshift=-1.5mm,yshift=1mm]$(AttrR\i)$);
				\coordinate (AttrR\i out) at ([xshift=-1.5mm,yshift=-1mm]$(AttrR\i)$);
				\coordinate (AttrR\i info) at ([xshift=1.5mm,yshift=0mm]$(AttrR\i)$);
			}}{}

\ifthenelse {\nAttrL > 0}
			{\foreach \i in {1,...,\nAttrL}{
				\coordinate (AttrL\i in) at ([xshift=1.5mm,yshift=1mm]$(AttrL\i)$);
				\coordinate (AttrL\i out) at ([xshift=1.5mm,yshift=-1mm]$(AttrL\i)$);
				\coordinate (AttrL\i info) at ([xshift=-1.5mm,yshift=0mm]$(AttrL\i)$);
			}}{}

\ifthenelse {\nMoveTo > 0}
			{\foreach \i in {1,...,\nMoveTo}{
				\coordinate (MoveTo\i inR) at ($(MoveTo\i .north west)!0.33!(MoveTo\i .north east)$);
				\coordinate (MoveTo\i inL) at ($(MoveTo\i .north west)!0.66!(MoveTo\i .north east)$);
				\coordinate (MoveTo\i outR) at ($(MoveTo\i .south west)!0.5!(MoveTo\i .south east)$);
			}}{}

\ifthenelse {\nPick > 0}
			{\foreach \i in {1,...,\nPick}{
					\coordinate (Pick\i inR) at ($(Pick\i .north west)!0.33!(Pick\i .north east)$);
					\coordinate (Pick\i inB) at ($(Pick\i .north west)!0.66!(Pick\i .north east)$);
					\coordinate (Pick\i outR) at ($(Pick\i .south west)!0.5!(Pick\i .south east)$);
			}}{}

\ifthenelse {\nPlace > 0}
			{\foreach \i in {1,...,\nPlace}{
					\coordinate (Place\i inR) at ($(Place\i .north west)!0.5!(Place\i .north east)$);
					\coordinate (Place\i outR) at ($(Place\i .south west)!0.33!(Place\i .south east)$);
					\coordinate (Place\i outB) at ($(Place\i .south west)!0.66!(Place\i .south east)$);
			}}{}
				
				% Draw edges
				\path[every node/.style={font=\scriptsize}]
				(in1) edge[ob] node[at start, yshift=\port] {robot} (MoveTo1inR)
				(in2) edge[info] node[at start, yshift=\port] {location} (Circ1)
				(Circ1) edge[info, out=180] node[] {} (MoveTo1inL)
				(Circ1) edge[info, out=0, in=0] node[] {} (AttrR1)
				(MoveTo1) edge[ob] node[] {} (AttrR1in)
				(AttrR1out) edge[ob] node[at end, yshift=-\port] {robot} (out1)
				;
				\end{tikzpicture}
			\end{tabular}
			&
			\begin{tabular}{c}
				=
			\end{tabular}
			&
			\begin{tabular}{c}
				\begin{tikzpicture}
				\tikzset{
	% Nodes
	AttrR/.style = {font={\Large$\rhd$}},
	AttrL/.style = {font={\Large$\lhd$}},
	MoveTo/.style = {draw,font={\MoveTo}},
	Pick/.style = {draw,font={\Pick}},
	Place/.style = {draw,font={\Place}},
	Circ/.style = {draw,circle,scale=0.75},
	Rect/.style = {draw,rectangle},
	%Edges
	ob/.style = {out=-90, in=90, looseness=1},
	info/.style = {dashed,out=-90, in=90, looseness=1},
	attrR/.style = {dashed,out=0, in=90, looseness=1},
	attrL/.style = {dashed,out=180, in=90, looseness=1},
	attrRL/.style = {dashed,out=0, in=180, looseness=1},
	attrLR/.style = {dashed,out=180, in=0, looseness=1},
%	attr/.style = {dashed,out=0, in=90, looseness=1},
%	urattr/.style = {dashed,out=-90, in=180, looseness=1},
%	ldattr/.style = {dashed,out=180, in=90, looseness=1},
%	ulattr/.style = {dashed,out=-90, in=0, looseness=1},
%	lrattr/.style = {dashed,out=180, in=0, looseness=1,shorten >=-8pt, shorten <=-8pt},
}

\pgfmathsetmacro{\port}{3mm}

				% Frame
				\coordinate (NW) at (0,2); \coordinate (NE) at (4,2); 
				\coordinate (SW) at (0,0); \coordinate (SE) at (4,0);
				%\draw[black] (NW) rectangle (SE);
				
				%Nodes
				\node[MoveTo] (MoveTo1) at (2, 1) {};
				
				% Autocreate port names
				\pgfmathsetmacro{\diagramIn}{2}
				\pgfmathsetmacro{\diagramOut}{1}
				\pgfmathsetmacro{\nAttrL}{0}
				\pgfmathsetmacro{\nAttrR}{1}
				\pgfmathsetmacro{\nMoveTo}{1}
				\pgfmathsetmacro{\nPick}{0}
				\pgfmathsetmacro{\nPlace}{0}
				\foreach \x [evaluate=\x as \dx using \x/(\diagramIn + 1)] in {1,...,\diagramIn}{
	\coordinate (in\x) at ([yshift=\port]$(NW)!\dx!(NE)$);
}

\foreach \x [evaluate=\x as \dx using \x/(\diagramOut + 1)] in {1,...,\diagramOut}{
	\coordinate (out\x) at ([yshift=-\port]$(SW)!\dx!(SE)$);
}

\ifthenelse	{\nAttrR > 0}
			{\foreach \i in {1,...,\nAttrR}{
				\coordinate (AttrR\i in) at ([xshift=-1.5mm,yshift=1mm]$(AttrR\i)$);
				\coordinate (AttrR\i out) at ([xshift=-1.5mm,yshift=-1mm]$(AttrR\i)$);
				\coordinate (AttrR\i info) at ([xshift=1.5mm,yshift=0mm]$(AttrR\i)$);
			}}{}

\ifthenelse {\nAttrL > 0}
			{\foreach \i in {1,...,\nAttrL}{
				\coordinate (AttrL\i in) at ([xshift=1.5mm,yshift=1mm]$(AttrL\i)$);
				\coordinate (AttrL\i out) at ([xshift=1.5mm,yshift=-1mm]$(AttrL\i)$);
				\coordinate (AttrL\i info) at ([xshift=-1.5mm,yshift=0mm]$(AttrL\i)$);
			}}{}

\ifthenelse {\nMoveTo > 0}
			{\foreach \i in {1,...,\nMoveTo}{
				\coordinate (MoveTo\i inR) at ($(MoveTo\i .north west)!0.33!(MoveTo\i .north east)$);
				\coordinate (MoveTo\i inL) at ($(MoveTo\i .north west)!0.66!(MoveTo\i .north east)$);
				\coordinate (MoveTo\i outR) at ($(MoveTo\i .south west)!0.5!(MoveTo\i .south east)$);
			}}{}

\ifthenelse {\nPick > 0}
			{\foreach \i in {1,...,\nPick}{
					\coordinate (Pick\i inR) at ($(Pick\i .north west)!0.33!(Pick\i .north east)$);
					\coordinate (Pick\i inB) at ($(Pick\i .north west)!0.66!(Pick\i .north east)$);
					\coordinate (Pick\i outR) at ($(Pick\i .south west)!0.5!(Pick\i .south east)$);
			}}{}

\ifthenelse {\nPlace > 0}
			{\foreach \i in {1,...,\nPlace}{
					\coordinate (Place\i inR) at ($(Place\i .north west)!0.5!(Place\i .north east)$);
					\coordinate (Place\i outR) at ($(Place\i .south west)!0.33!(Place\i .south east)$);
					\coordinate (Place\i outB) at ($(Place\i .south west)!0.66!(Place\i .south east)$);
			}}{}
				
				% Draw edges
				\path[every node/.style={font=\scriptsize}]
				(in1) edge[ob] node[at start, yshift=\port] {robot} (MoveTo1inR)
				(in2) edge[info] node[at start, yshift=\port] {location} (MoveTo1inL)
				(MoveTo1) edge[ob] node[at end, yshift=-\port] {robot} (out1)
				;
				\end{tikzpicture}
			\end{tabular}
	\end{tabular}\label{eq:moveto}\end{equation}
Here (and throughout the paper) solid lines represent entities and dashed lines represent data. The equation defines a post-condition of the operation: after executing the process $\MoveTo$, the robot has moved to the specified location. Given similar axioms for $\Pick$ and $\Place$ operations, we can use string diagrams to prove that a sequence of operations is well-defined or validates a desired specification, as shown in Figure \ref{fig:proofs}.
These specific examples are formalized in Subsection \ref{sec:example}.

	\begin{figure}
		\centering
	\begin{tabular}{|c|c|}\hline
	\scalebox{.66}{
	\begin{tabular}{ccc}
	\begin{tabular}{c}\begin{tikzpicture}
		\input{mystyles}
		% Frame
		\coordinate (NW) at (0,4); \coordinate (NE) at (8,4); 
		\coordinate (SW) at (0,0); \coordinate (SE) at (8,0);
		%\draw[black] (NW) rectangle (SE);
		
		%Nodes
		\node[AttrL] (AttrL1) at ({$(NW)!0.63!(NE)$} |- {$(SW)!0.85!(NW)$}) {};
		\node[MoveTo] (MoveTo1) at ({$(NW)!0.4!(NE)$} |- {$(SW)!0.6!(NW)$}) {};
		\node[Pick] (Pick1) at ({$(NW)!0.5!(NE)$} |- {$(SW)!0.2!(NW)$}) {};
		
		% Autocreate port names
		\pgfmathsetmacro{\diagramIn}{2}
		\pgfmathsetmacro{\diagramOut}{1}	
		\pgfmathsetmacro{\nAttrL}{1}
		\pgfmathsetmacro{\nAttrR}{0}
		\pgfmathsetmacro{\nMoveTo}{1}
		\pgfmathsetmacro{\nPick}{1}
		\pgfmathsetmacro{\nPlace}{0}
		\input{pgfpositions}
		
		% Draw edges
		\path[every node/.style={font=\scriptsize}]
		(in1) edge[ob] node[at start,yshift=\port] {robot} (MoveTo1inR)
		(in2) edge[ob] node[at start,yshift=\port] {ball} (AttrL1in)
		(AttrL1info) edge[attrL] (MoveTo1inL)
		(AttrL1out) edge[ob] (Pick1inB)
		(MoveTo1outR) edge[ob] (Pick1inR)
		(Pick1outR) edge[ob] node[at end,yshift=-\port] {robot} (out1)
		;
		
		\end{tikzpicture}\end{tabular}
	&
	\begin{tabular}{c}
		\LARGE =
	\end{tabular}
	&
	\begin{tabular}{c}	\begin{tikzpicture}
		\input{mystyles}
		% Frame
		\coordinate (NW) at (0,4); \coordinate (NE) at (6,4); 
		\coordinate (SW) at (0,0); \coordinate (SE) at (6,0);
		%\draw[black] (NW) rectangle (SE);
		
		%Nodes
		\node[AttrL] (AttrL1) at ({$(NW)!0.72!(NE)$} |- {$(SW)!0.85!(NW)$}) {};
		\node[Circ] (Circ1) at ({$(NW)!0.5!(NE)$} |- {$(SW)!0.85!(NW)$}) {};
		\node[AttrR] (AttrR1) at ({$(NW)!0.34!(NE)$} |- {$(SW)!0.4!(NW)$}) {};
		\node[AttrL] (AttrL2) at ({$(NW)!0.66!(NE)$} |- {$(SW)!0.4!(NW)$}) {};
		\node[Circ] (Circ2) at ({$(NW)!0.5!(NE)$} |- {$(SW)!0.4!(NW)$}) {};
		\node[MoveTo] (MoveTo1) at ({$(NW)!0.3!(NE)$} |- {$(SW)!0.6!(NW)$}) {};
		\node[Pick] (Pick1) at ({$(NW)!0.5!(NE)$} |- {$(SW)!0.1!(NW)$}) {};
		
		% Autocreate port names
		\pgfmathsetmacro{\diagramIn}{3}
		\pgfmathsetmacro{\diagramOut}{1}	
		\pgfmathsetmacro{\nAttrL}{2}
		\pgfmathsetmacro{\nAttrR}{1}
		\pgfmathsetmacro{\nMoveTo}{1}
		\pgfmathsetmacro{\nPick}{1}
		\pgfmathsetmacro{\nPlace}{0}
		\input{pgfpositions}
		
		% Draw edges
		\path[every node/.style={font=\scriptsize}]
		(in1) edge[ob] node[at start,yshift=\port] {robot} (MoveTo1inR)
		(in3) edge[ob] node[at start,yshift=\port] {ball} (AttrL1in)
		(AttrL1info) edge[attrLR] (Circ1)
		(Circ1) edge[attrL] (MoveTo1inL)
		(Circ1) edge[info] (Circ2)
		(AttrL1out) edge[ob] (AttrL2in)
		(MoveTo1outR) edge[ob] (AttrR1in)
		(AttrR1info) edge[attrRL] (Circ2)
		(AttrL2info) edge[attrLR] (Circ2)
		(AttrR1out) edge[ob] (Pick1inR)
		(AttrL2out) edge[ob] (Pick1inB)
		(Pick1outR) edge[ob] node[at end,yshift=-\port] {robot} (out1)
		;
		
		\end{tikzpicture}\end{tabular}
	\end{tabular}
	}
	&
	\scalebox{.66}{
	\begin{tabular}{ccc}
	\begin{tabular}{c}\begin{tikzpicture}
		\input{mystyles}
		% Frame
		\coordinate (NW) at (0,4); \coordinate (NE) at (6,4); 
		\coordinate (SW) at (0,0); \coordinate (SE) at (6,0);
		%\draw[black] (NW) rectangle (SE);
		
		%Nodes
		\node[MoveTo] (MoveTo1) at ({$(NW)!0.5!(NE)$} |- {$(SW)!0.7!(NW)$}) {};
		\node[Place] (Place1) at ({$(NW)!0.5!(NE)$} |- {$(SW)!0.3!(NW)$}) {};
		
		% Autocreate port names
		\pgfmathsetmacro{\diagramIn}{2}
		\pgfmathsetmacro{\diagramOut}{2}	
		\pgfmathsetmacro{\nAttrL}{0}
		\pgfmathsetmacro{\nAttrR}{0}
		\pgfmathsetmacro{\nMoveTo}{1}
		\pgfmathsetmacro{\nPick}{1}
		\pgfmathsetmacro{\nPlace}{1}
		\input{pgfpositions}
		
		% Draw edges
		\path[every node/.style={font=\scriptsize}]
		(in1) edge[ob] node[at start,yshift=\port] {robot} (MoveTo1inR)
		(in2) edge[info] node[at start,yshift=\port] {goal} (MoveTo1inL)
		(MoveTo1outR) edge[ob] (Place1inR)
		(Place1outR) edge[ob] node[at end,yshift=-\port] {robot} (out1)
		(Place1outB) edge[ob] node[at end,yshift=-\port] {ball} (out2)
		;
		
		\end{tikzpicture}\end{tabular}
	&
	\begin{tabular}{c}
		\LARGE =
	\end{tabular}
	&
	\begin{tabular}{c}\begin{tikzpicture}
		\input{mystyles}
		% Frame
		\coordinate (NW) at (0,4); \coordinate (NE) at (6,4); 
		\coordinate (SW) at (0,0); \coordinate (SE) at (6,0);
		%\draw[black] (NW) rectangle (SE);
		
		%Nodes
		\node[Circ] (Circ1) at ({$(NW)!0.66!(NE)$} |- {$(SW)!0.9!(NW)$}) {};
		\node[MoveTo] (MoveTo1) at ({$(NW)!0.45!(NE)$} |- {$(SW)!0.7!(NW)$}) {};
		\node[AttrR] (AttrR1) at ({$(NW)!0.48!(NE)$} |- {$(SW)!0.5!(NW)$}) {};
		\node[Circ] (Circ2) at ({$(NW)!0.66!(NE)$} |- {$(SW)!0.5!(NW)$}) {};
		\node[Place] (Place1) at ({$(NW)!0.45!(NE)$} |- {$(SW)!0.3!(NW)$}) {};
		\node[AttrR] (AttrR2) at ({$(NW)!0.58!(NE)$} |- {$(SW)!0.1!(NW)$}) {};
		
		% Autocreate port names
		\pgfmathsetmacro{\diagramIn}{2}
		\pgfmathsetmacro{\diagramOut}{2}	
		\pgfmathsetmacro{\nAttrL}{1}
		\pgfmathsetmacro{\nAttrR}{2}
		\pgfmathsetmacro{\nMoveTo}{1}
		\pgfmathsetmacro{\nPick}{1}
		\pgfmathsetmacro{\nPlace}{1}
		\input{pgfpositions}
		
		% Draw edges
		\path[every node/.style={font=\scriptsize}]
		(in1) edge[ob] node[at start,yshift=\port] {robot} (MoveTo1inR)
		(in2) edge[info] node[at start,yshift=\port] {goal} (Circ1)
		(Circ1) edge[attrL] (MoveTo1inL)
		(Circ1) edge[info] (Circ2)
		(MoveTo1outR) edge[ob] (AttrR1in)
		(AttrR1out) edge[ob] (Place1inR)
		(AttrR1info) edge[attrR,in=180] (Circ2)
		(AttrR2info) edge[attrR,in=-90] (Circ2)
		(Place1outR) edge[ob] node[at end,yshift=-\port] {robot} (out1)
		(Place1outB) edge[ob] (AttrR2in)
		(AttrR2out) edge[ob] node[at end,yshift=-\port] {ball} (out2)
		;
		
		\end{tikzpicture}\end{tabular}
	\end{tabular}
	}\\
	(a) & (b)\\\hline
	\end{tabular}
	\caption[caption]{(a) This equation guarantees that the preconditions of $\Pick$ are satisfied. \\\hspace{\textwidth} (b) This equation guarantees the desired specification (ball at goal).}
	\label{fig:proofs}
	\end{figure}

A few questions arise when trying to formalize this.
First, how are informational resources supposed to behave in general?
To model such resources, we borrow the notion of a \emph{data service} from \cite{pavlovic2013monoidal}; we summarize this in Definition \ref{def:dataservices}.
Data services are algebraic structures (defined internally to SMCs) which formalize the operations of \emph{filtering}, \emph{copying}, and \emph{deleting} pieces of data.

Furthermore, what does it mean for an entity to have a datum as an attribute (e.g.\ in the way that the robot above ``has a location'')?
We argue that this question can be answered by requiring that the data service associated to the informational object ``acts on'' the other object in a suitable sense. 
To provide more control over the behavior of these attributes, we define the notion of an \emph{attribute structure} on a SMC.
These concepts are discussed in Subsection \ref{sec:generalcase}.

When considering the semantics of these categories, it is often the case that certain morphisms turn out to be ``partially defined.'' 
This manifests itself through a partial order structure on hom-sets, or more precisely an enrichment over $\Poset$, the category of partially-ordered sets (posets) and monotone increasing functions.
Due to similarities between these partial orders across semantics, we find it useful and illuminating to enrich the \emph{syntax} categories over $\Poset$ as well. 
By doing so we are able to impose extra conditions on these attribute structures that mimic the ways users expect informational resources to behave.
The details of this enrichment are spelled out in Subsection \ref{sec:posetcase}.

In our motivating example, we use these categories with attributes and their diagrams to coordinate two semantic models at different levels of abstraction, based on the principle of functorial semantics. First, we define a Boolean semantics based on the Planning Domain Definition Language (PDDL) \cite{PDDL}; a presentation of the categorical syntax can be used to generate a PDDL domain, and the solution to a PDDL problem can be used to define an associated string diagram.
Once a high-level plan has been identified, a second mapping to a category of geometric semantics can be used to parameterize a path planning algorithm. We give sketches of these categories and mappings in Section \ref{sec:semantics}.

We seek to achieve two principal goals in this paper.
First, we aim to develop rigorous foundations for the forthcoming paper \cite{robotspaper}, which explores similar ideas with a greater focus on the engineering than on the mathematics.
Second, we seek to present interesting examples of categorical modeling in engineering.
It is our hope that the applications presented in this paper will motivate other researchers to investigate the connections between category theory and engineering models.

\section{Categories with Attributes}

Our goal in this section is to develop a definition of a \emph{category with attributes} -- that is, a (symmetric monoidal) category in which objects can have ``information'' (defined internally to the category) attached to them in some way.
We develop one notion of a category with attributes in Subsection \ref{sec:generalcase}.
In Subsection \ref{sec:posetcase}, we adapt this notion to the case where the categories involved are $\Poset$-enriched.
This additional structure allows us to impose conditions on the relationships between certain morphisms, allowing us to connect categories with attributes to our intuition about partially-defined morphisms (e.g.\ partial functions). 

\subsection{General Case} \label{sec:generalcase}

For a formal categorical notion of information, we borrow the definition of a \emph{data service} from \cite{pavlovic2013monoidal}.\footnote{
Readers familiar with Frobenius algebras will observe that a data service is the same as a special commutative Frobenius algebra except for the fact that data services are not required to have units.
Non-unitality can be unavoidable in certain applications; for example, every object in $\mathsf{PartFn}$, the category of sets and partial functions, has a canonical data service structure, but this structure typically does not admit a unit.
}
This definition is restated as follows.
As hinted above, we will depict the underlying objects of data services with \emph{dashed lines} to contrast them with other objects.

% =======================================================================================================================
% SHRINK ME WHEN NAVIGATING THE PAPER
% =======================================================================================================================
\begin{defn}[\cite{pavlovic2013monoidal}] \label{def:dataservices}
A \emph{data service} in a SMC $(\CC, \otimes, I)$ consists of
\begin{itemize}
    \item An object $D$ in $\CC$;
    \item A \emph{multiplication morphism} $\mu: D \otimes D \to D$;
    \item A \emph{comultiplication morphism} $\delta: D \to D \otimes D$; and
    \item A \emph{counit morphism} $\epsilon: D \to I$.
\end{itemize}
The morphisms of a data service are depicted as follows:
\begin{equation}
\end{equation}
These morphisms are also required to satisfy several axioms:
\begin{itemize}
    \item $(V, \mu)$ is a commutative semigroup object in $\CC$;
    \item $(V, \delta, \epsilon)$ is a commutative comonoid object in $\CC$;
    \item $\mu$ and $\delta$ satisfy the \emph{Frobenius laws}:
\begin{equation}%
\end{tabular}\end{equation}
\end{itemize}
When we have no reason to explicitly reference the morphisms $\mu$, $\delta$, and $\epsilon$, we will refer to the data service $(D, \mu, \delta, \epsilon)$ simply as $D$.
\end{defn}

% THIS CAN BE DROPPED / CONDENSED AS NECESSARY

Some discussion of the intuition behind data services is in order.
As mentioned above, an object with a data service structure can be thought of as a value or datum.
From this perspective, the comultiplication morphism can be viewed as creating a copy of the input value (where the copy is indistinguishable from the original), while the counit morphism can be viewed as deleting or forgetting about the input value.
Meanwhile, the multiplication morphism can be thought of as filtering for equality: given two values, the multiplication morphism checks whether or not they're the same and returns the common value of both (if they're the same) or nothing (if they aren't).
The Frobenius, special, semigroup, and comonoid axioms formalize this intuition.

For any SMC $\CC$, we can define a category $\Data(\CC)$ of data services in $\CC$.
The ``obvious'' definition of a morphism of data services $D \to D'$ (i.e.\ a morphism $D \to D'$ in $\CC$ that is both a semigroup homomorphism and a comonoid homomorphism) will turn out to be too restrictive for our eventual goals.
In this general case, we are unable to find any nontrivial useful ways to weaken this ``obvious'' definition, so we opt for the trivial solution, allowing any morphism between the underlying objects of two data services to be considered a ``morphism of data services.''
(Contrast this with the discussion in Subsection \ref{sec:posetcase}, where we are able to formulate acceptable conditions).
This is recorded in the following definition.

\begin{defn} \label{def:dataC}
Let $\CC$ be a SMC. The \emph{category of data services} in $\CC$, denoted $\Data(\CC)$, is the category where:
\begin{itemize}
    \item Objects are data services in $\CC$
    \item $\Data(\CC)(D, D') = \CC(D, D')$, i.e.\ morphisms in $\Data(\CC)$ are arbitrary morphisms between the underlying objects of the domain / codomain.
\end{itemize}
There exists an obvious (fully faithful) forgetful functor $U: \Data(\CC) \to \CC$.
\end{defn}

In order to define categories with attributes, we still need to determine a proper definition for ``attributes'' themselves.
To an initial approximation, we consider an attribute of an object $M$ to be a data service $D$ associated to $D$ in some way.
Given an instance of $M$, we should be able to retrieve the associated instance of $D$.
This notion can be formalized by the morphism $\gamma: M \to M \otimes D$ depicted in \eqref{eq:gammaphi}.\footnote{
An attribute is not determined by the pair $(M, D)$, so our depiction here is lacking insofar as it does not represent the specific choice of attribute for $(M, D)$.
Despite this shortcoming, the depiction here is still useful for brevity, and we have no need to explicitly consider different attributes with the same $D$ and $M$.
}
Furthermore, we should be able to ``filter attributes for equality,'' checking whether or not the instance of $D$ associated to the instance of $M$ agrees with some arbitrary instance of $D$.
This is formalized by the morphism $\phi: M \otimes D \to M$ depicted in \eqref{eq:gammaphi}.

% =======================================================================================================================
% SHRINK ME WHEN NAVIGATING THE PAPER
% =======================================================================================================================
\begin{equation} \label{eq:gammaphi}
\begin{tabular}{cccc} 
\begin{tabular}{c}
$\gamma := $
\end{tabular}
&
\begin{tabular}{c}
\begin{tikzpicture}
    \tikzset{
	% Nodes
	AttrR/.style = {font={\Large$\rhd$}},
	AttrL/.style = {font={\Large$\lhd$}},
	MoveTo/.style = {draw,font={\MoveTo}},
	Pick/.style = {draw,font={\Pick}},
	Place/.style = {draw,font={\Place}},
	Circ/.style = {draw,circle,scale=0.75},
	Rect/.style = {draw,rectangle},
	%Edges
	ob/.style = {out=-90, in=90, looseness=1},
	info/.style = {dashed,out=-90, in=90, looseness=1},
	attrR/.style = {dashed,out=0, in=90, looseness=1},
	attrL/.style = {dashed,out=180, in=90, looseness=1},
	attrRL/.style = {dashed,out=0, in=180, looseness=1},
	attrLR/.style = {dashed,out=180, in=0, looseness=1},
%	attr/.style = {dashed,out=0, in=90, looseness=1},
%	urattr/.style = {dashed,out=-90, in=180, looseness=1},
%	ldattr/.style = {dashed,out=180, in=90, looseness=1},
%	ulattr/.style = {dashed,out=-90, in=0, looseness=1},
%	lrattr/.style = {dashed,out=180, in=0, looseness=1,shorten >=-8pt, shorten <=-8pt},
}

\pgfmathsetmacro{\port}{3mm}

	% Frame
	\coordinate (NW) at (0,1); \coordinate (NE) at (3,1); 
	\coordinate (SW) at (0,0); \coordinate (SE) at (3,0);
	%\draw[black] (NW) rectangle (SE);
	
	%Nodes
	\node[AttrR] (AttrR1) at (1.15, .5) {};
	
	% Autocreate port names
	\pgfmathsetmacro{\diagramIn}{2}
	\pgfmathsetmacro{\diagramOut}{2}	
	\pgfmathsetmacro{\nAttrL}{0}
	\pgfmathsetmacro{\nAttrR}{1}
	\pgfmathsetmacro{\nMoveTo}{0}
	\pgfmathsetmacro{\nPick}{0}
	\pgfmathsetmacro{\nPlace}{0}
	\foreach \x [evaluate=\x as \dx using \x/(\diagramIn + 1)] in {1,...,\diagramIn}{
	\coordinate (in\x) at ([yshift=\port]$(NW)!\dx!(NE)$);
}

\foreach \x [evaluate=\x as \dx using \x/(\diagramOut + 1)] in {1,...,\diagramOut}{
	\coordinate (out\x) at ([yshift=-\port]$(SW)!\dx!(SE)$);
}

\ifthenelse	{\nAttrR > 0}
			{\foreach \i in {1,...,\nAttrR}{
				\coordinate (AttrR\i in) at ([xshift=-1.5mm,yshift=1mm]$(AttrR\i)$);
				\coordinate (AttrR\i out) at ([xshift=-1.5mm,yshift=-1mm]$(AttrR\i)$);
				\coordinate (AttrR\i info) at ([xshift=1.5mm,yshift=0mm]$(AttrR\i)$);
			}}{}

\ifthenelse {\nAttrL > 0}
			{\foreach \i in {1,...,\nAttrL}{
				\coordinate (AttrL\i in) at ([xshift=1.5mm,yshift=1mm]$(AttrL\i)$);
				\coordinate (AttrL\i out) at ([xshift=1.5mm,yshift=-1mm]$(AttrL\i)$);
				\coordinate (AttrL\i info) at ([xshift=-1.5mm,yshift=0mm]$(AttrL\i)$);
			}}{}

\ifthenelse {\nMoveTo > 0}
			{\foreach \i in {1,...,\nMoveTo}{
				\coordinate (MoveTo\i inR) at ($(MoveTo\i .north west)!0.33!(MoveTo\i .north east)$);
				\coordinate (MoveTo\i inL) at ($(MoveTo\i .north west)!0.66!(MoveTo\i .north east)$);
				\coordinate (MoveTo\i outR) at ($(MoveTo\i .south west)!0.5!(MoveTo\i .south east)$);
			}}{}

\ifthenelse {\nPick > 0}
			{\foreach \i in {1,...,\nPick}{
					\coordinate (Pick\i inR) at ($(Pick\i .north west)!0.33!(Pick\i .north east)$);
					\coordinate (Pick\i inB) at ($(Pick\i .north west)!0.66!(Pick\i .north east)$);
					\coordinate (Pick\i outR) at ($(Pick\i .south west)!0.5!(Pick\i .south east)$);
			}}{}

\ifthenelse {\nPlace > 0}
			{\foreach \i in {1,...,\nPlace}{
					\coordinate (Place\i inR) at ($(Place\i .north west)!0.5!(Place\i .north east)$);
					\coordinate (Place\i outR) at ($(Place\i .south west)!0.33!(Place\i .south east)$);
					\coordinate (Place\i outB) at ($(Place\i .south west)!0.66!(Place\i .south east)$);
			}}{}
	
	% Draw edges
	\path[every node/.style={font=\scriptsize}]
	(in1) edge[ob] node[at start, yshift=\port] {$M$} (AttrR1in)
	(AttrR1out) edge[ob] node[at end, yshift=-\port] {$M$} (out1)
	(AttrR1) edge[attrR] node[at end, yshift=-\port] {$D$} (out2)
	;
\end{tikzpicture}
\end{tabular};
&
\begin{tabular}{c}
$\phi := $
\end{tabular}
&
\begin{tabular}{c}
    \begin{tikzpicture}
        \tikzset{
	% Nodes
	AttrR/.style = {font={\Large$\rhd$}},
	AttrL/.style = {font={\Large$\lhd$}},
	MoveTo/.style = {draw,font={\MoveTo}},
	Pick/.style = {draw,font={\Pick}},
	Place/.style = {draw,font={\Place}},
	Circ/.style = {draw,circle,scale=0.75},
	Rect/.style = {draw,rectangle},
	%Edges
	ob/.style = {out=-90, in=90, looseness=1},
	info/.style = {dashed,out=-90, in=90, looseness=1},
	attrR/.style = {dashed,out=0, in=90, looseness=1},
	attrL/.style = {dashed,out=180, in=90, looseness=1},
	attrRL/.style = {dashed,out=0, in=180, looseness=1},
	attrLR/.style = {dashed,out=180, in=0, looseness=1},
%	attr/.style = {dashed,out=0, in=90, looseness=1},
%	urattr/.style = {dashed,out=-90, in=180, looseness=1},
%	ldattr/.style = {dashed,out=180, in=90, looseness=1},
%	ulattr/.style = {dashed,out=-90, in=0, looseness=1},
%	lrattr/.style = {dashed,out=180, in=0, looseness=1,shorten >=-8pt, shorten <=-8pt},
}

\pgfmathsetmacro{\port}{3mm}

    	% Frame
    	\coordinate (NW) at (0,1); \coordinate (NE) at (3,1); 
    	\coordinate (SW) at (0,0); \coordinate (SE) at (3,0);
    	%\draw[black] (NW) rectangle (SE);
    	
    	%Nodes
    	\node[AttrR] (AttrR1) at (1.15, 0.5) {};
    	
    	% Autocreate port names
    	\pgfmathsetmacro{\diagramIn}{2}
    	\pgfmathsetmacro{\diagramOut}{2}	
    	\pgfmathsetmacro{\nAttrL}{0}
    	\pgfmathsetmacro{\nAttrR}{1}
    	\pgfmathsetmacro{\nMoveTo}{0}
    	\pgfmathsetmacro{\nPick}{0}
    	\pgfmathsetmacro{\nPlace}{0}
    	\foreach \x [evaluate=\x as \dx using \x/(\diagramIn + 1)] in {1,...,\diagramIn}{
	\coordinate (in\x) at ([yshift=\port]$(NW)!\dx!(NE)$);
}

\foreach \x [evaluate=\x as \dx using \x/(\diagramOut + 1)] in {1,...,\diagramOut}{
	\coordinate (out\x) at ([yshift=-\port]$(SW)!\dx!(SE)$);
}

\ifthenelse	{\nAttrR > 0}
			{\foreach \i in {1,...,\nAttrR}{
				\coordinate (AttrR\i in) at ([xshift=-1.5mm,yshift=1mm]$(AttrR\i)$);
				\coordinate (AttrR\i out) at ([xshift=-1.5mm,yshift=-1mm]$(AttrR\i)$);
				\coordinate (AttrR\i info) at ([xshift=1.5mm,yshift=0mm]$(AttrR\i)$);
			}}{}

\ifthenelse {\nAttrL > 0}
			{\foreach \i in {1,...,\nAttrL}{
				\coordinate (AttrL\i in) at ([xshift=1.5mm,yshift=1mm]$(AttrL\i)$);
				\coordinate (AttrL\i out) at ([xshift=1.5mm,yshift=-1mm]$(AttrL\i)$);
				\coordinate (AttrL\i info) at ([xshift=-1.5mm,yshift=0mm]$(AttrL\i)$);
			}}{}

\ifthenelse {\nMoveTo > 0}
			{\foreach \i in {1,...,\nMoveTo}{
				\coordinate (MoveTo\i inR) at ($(MoveTo\i .north west)!0.33!(MoveTo\i .north east)$);
				\coordinate (MoveTo\i inL) at ($(MoveTo\i .north west)!0.66!(MoveTo\i .north east)$);
				\coordinate (MoveTo\i outR) at ($(MoveTo\i .south west)!0.5!(MoveTo\i .south east)$);
			}}{}

\ifthenelse {\nPick > 0}
			{\foreach \i in {1,...,\nPick}{
					\coordinate (Pick\i inR) at ($(Pick\i .north west)!0.33!(Pick\i .north east)$);
					\coordinate (Pick\i inB) at ($(Pick\i .north west)!0.66!(Pick\i .north east)$);
					\coordinate (Pick\i outR) at ($(Pick\i .south west)!0.5!(Pick\i .south east)$);
			}}{}

\ifthenelse {\nPlace > 0}
			{\foreach \i in {1,...,\nPlace}{
					\coordinate (Place\i inR) at ($(Place\i .north west)!0.5!(Place\i .north east)$);
					\coordinate (Place\i outR) at ($(Place\i .south west)!0.33!(Place\i .south east)$);
					\coordinate (Place\i outB) at ($(Place\i .south west)!0.66!(Place\i .south east)$);
			}}{}
    	
    	% Draw edges
    	\path[every node/.style={font=\scriptsize}]
    	(in1) edge[ob] node[at start, yshift=\port] {$M$} (AttrR1in)
    	(in2) edge[info, in=0] node[at start, yshift=\port] {$D$} (AttrR1)
    	(AttrR1out) edge[ob] node[at end, yshift=-\port] {$M$} (out1)
    	;
    \end{tikzpicture}
\end{tabular}
\end{tabular}\end{equation}

These morphisms should be expected to satisfy some compatibility conditions with the data service structures on the object $D$.
For example, we expect that retrieving a datum twice in a row should be the same as retrieving the datum once and then copying the result.
A collection of similar axioms are presented in the following definition.

% =======================================================================================================================
% SHRINK ME WHEN NAVIGATING THE PAPER
% =======================================================================================================================
\begin{defn} \label{def:dataserviceaction}
Let $\CC$ be a SMC, let $(D, \mu, \delta, \epsilon)$ be a data service in $\CC$, and let $M$ be any object of $\CC$.
A (right) \emph{data service action} of $D$ on $M$ is a pair $(\gamma, \phi)$ of morphisms $\gamma: M \to M \otimes D$ and $\phi: M \otimes D \to M$, depicted as in \eqref{eq:gammaphi}, such that
\begin{itemize}
    \item $\phi$ gives a (right) action of the semigroup object $(V, \mu)$ on $E$:
\begin{equation}
\end{equation}
    \item $\gamma$ gives a (right) action of the comonoid object $(M, \delta, \epsilon)$ on $E$:
\begin{equation}%
\end{tabular}\end{equation}   
\end{itemize}
\end{defn}

Remarkably, a data service action $(\gamma, \phi)$ is entirely determined by the comonoid action $\gamma$, as the following proposition shows.\footnote{
This proposition is closely related to results obtained in \cite{abrams1999modules} for Frobenius algebras.
Note that the asymmetry of the data service axioms prevents us from recovering $\gamma$ from $\phi$ alone in a data service, as in \cite{abrams1999modules}.
}

% =======================================================================================================================
% SHRINK ME WHEN NAVIGATING THE PAPER
% =======================================================================================================================
\begin{prop} \label{prop:comonoid}
Let $(D, \mu, \delta, \epsilon)$ be a data service and let $\gamma: M \to M \otimes D$ (depicted as above) be a (right) action of the comonoid object $(D, \delta, \epsilon)$ on the object $M$.
Then there exists a unique morphism $\phi: M \otimes D \to M$ such that $(\gamma, \phi)$ forms a data service action of $D$ on $M$.
This $\phi$ is defined as in \eqref{eq:phidef}.
\begin{equation} \label{eq:phidef}
\begin{tabular}{cc}
\begin{tabular}{c}
$\phi :=$
\end{tabular}
&
\begin{tabular}{c}
    \begin{tikzpicture}
        \tikzset{
	% Nodes
	AttrR/.style = {font={\Large$\rhd$}},
	AttrL/.style = {font={\Large$\lhd$}},
	MoveTo/.style = {draw,font={\MoveTo}},
	Pick/.style = {draw,font={\Pick}},
	Place/.style = {draw,font={\Place}},
	Circ/.style = {draw,circle,scale=0.75},
	Rect/.style = {draw,rectangle},
	%Edges
	ob/.style = {out=-90, in=90, looseness=1},
	info/.style = {dashed,out=-90, in=90, looseness=1},
	attrR/.style = {dashed,out=0, in=90, looseness=1},
	attrL/.style = {dashed,out=180, in=90, looseness=1},
	attrRL/.style = {dashed,out=0, in=180, looseness=1},
	attrLR/.style = {dashed,out=180, in=0, looseness=1},
%	attr/.style = {dashed,out=0, in=90, looseness=1},
%	urattr/.style = {dashed,out=-90, in=180, looseness=1},
%	ldattr/.style = {dashed,out=180, in=90, looseness=1},
%	ulattr/.style = {dashed,out=-90, in=0, looseness=1},
%	lrattr/.style = {dashed,out=180, in=0, looseness=1,shorten >=-8pt, shorten <=-8pt},
}

\pgfmathsetmacro{\port}{3mm}

    	% Frame
    	\coordinate (NW) at (0,1); \coordinate (NE) at (3,1); 
    	\coordinate (SW) at (0,0); \coordinate (SE) at (3,0);
    	%\draw[black] (NW) rectangle (SE);
    	
    	%Nodes
    	\node[AttrR] (AttrR1) at (0.75, 0.9) {};
    	\node[Circ] (Circ1) at (2, 0.7) {};
    	\node[Circ] (Circ2) at (2, 0.1) {};
    	
    	% Autocreate port names
    	\pgfmathsetmacro{\diagramIn}{4}
    	\pgfmathsetmacro{\diagramOut}{4}	
    	\pgfmathsetmacro{\nAttrL}{0}
    	\pgfmathsetmacro{\nAttrR}{1}
    	\pgfmathsetmacro{\nMoveTo}{0}
    	\pgfmathsetmacro{\nPick}{0}
    	\pgfmathsetmacro{\nPlace}{0}
    	\foreach \x [evaluate=\x as \dx using \x/(\diagramIn + 1)] in {1,...,\diagramIn}{
	\coordinate (in\x) at ([yshift=\port]$(NW)!\dx!(NE)$);
}

\foreach \x [evaluate=\x as \dx using \x/(\diagramOut + 1)] in {1,...,\diagramOut}{
	\coordinate (out\x) at ([yshift=-\port]$(SW)!\dx!(SE)$);
}

\ifthenelse	{\nAttrR > 0}
			{\foreach \i in {1,...,\nAttrR}{
				\coordinate (AttrR\i in) at ([xshift=-1.5mm,yshift=1mm]$(AttrR\i)$);
				\coordinate (AttrR\i out) at ([xshift=-1.5mm,yshift=-1mm]$(AttrR\i)$);
				\coordinate (AttrR\i info) at ([xshift=1.5mm,yshift=0mm]$(AttrR\i)$);
			}}{}

\ifthenelse {\nAttrL > 0}
			{\foreach \i in {1,...,\nAttrL}{
				\coordinate (AttrL\i in) at ([xshift=1.5mm,yshift=1mm]$(AttrL\i)$);
				\coordinate (AttrL\i out) at ([xshift=1.5mm,yshift=-1mm]$(AttrL\i)$);
				\coordinate (AttrL\i info) at ([xshift=-1.5mm,yshift=0mm]$(AttrL\i)$);
			}}{}

\ifthenelse {\nMoveTo > 0}
			{\foreach \i in {1,...,\nMoveTo}{
				\coordinate (MoveTo\i inR) at ($(MoveTo\i .north west)!0.33!(MoveTo\i .north east)$);
				\coordinate (MoveTo\i inL) at ($(MoveTo\i .north west)!0.66!(MoveTo\i .north east)$);
				\coordinate (MoveTo\i outR) at ($(MoveTo\i .south west)!0.5!(MoveTo\i .south east)$);
			}}{}

\ifthenelse {\nPick > 0}
			{\foreach \i in {1,...,\nPick}{
					\coordinate (Pick\i inR) at ($(Pick\i .north west)!0.33!(Pick\i .north east)$);
					\coordinate (Pick\i inB) at ($(Pick\i .north west)!0.66!(Pick\i .north east)$);
					\coordinate (Pick\i outR) at ($(Pick\i .south west)!0.5!(Pick\i .south east)$);
			}}{}

\ifthenelse {\nPlace > 0}
			{\foreach \i in {1,...,\nPlace}{
					\coordinate (Place\i inR) at ($(Place\i .north west)!0.5!(Place\i .north east)$);
					\coordinate (Place\i outR) at ($(Place\i .south west)!0.33!(Place\i .south east)$);
					\coordinate (Place\i outB) at ($(Place\i .south west)!0.66!(Place\i .south east)$);
			}}{}
    	
    	% Draw edges
    	\path[every node/.style={font=\scriptsize}]
    	(in1) edge[ob] node[at start, yshift=\port] {$M$} (AttrR1in)
    	(in4) edge[info, in=0] node[at start, yshift=\port] {$D$} (Circ1)
    	(AttrR1) edge[info, out=0, in=180] node[yshift=\port] {$D$} (Circ1)
    	(Circ1) edge[info] node[right] {$D$} (Circ2)
    	(AttrR1out) edge[ob] node[at end, yshift=-\port] {$M$} (out1)
    	;
    \end{tikzpicture}
\end{tabular}
\end{tabular}\end{equation}
\end{prop}

% CONDENSE THE FOLLOWING AS NECESSARY

The proof that such a $\phi$ is unique is fairly straightforward; in fact, one only needs to use the equality between the first and third terms of \eqref{eq:frobeniusaction} and the fact that $D$ is a data service.
A series of immediate diagram chases then establishes all of the other laws of Definition \ref{def:dataserviceaction} (provided that $\gamma$ gives a comonoid action of $(D, \delta, \epsilon)$ on $M$).
In particular, the data service action axioms are overdetermined, so that some of the less immediately intuitive laws (e.g.\ the equality between the first and second terms of \eqref{eq:frobeniusaction}) follows from the other assumptions

Given a data service action $(\gamma, \phi)$ of some data service $D$ on some object $M$, we can define several other interesting morphisms.
For example, because $\CC$ is symmetric and $D$ is a commutative semigroup / comonoid object, we obtain a pair of morphisms $(\gamma', \phi')$ giving a ``left data service action'' of $D$ on $M$, defined by (post- or pre-)composing $\gamma$ and $\phi$ with the appropriate braiding.
We depict these $\gamma'$ and $\phi'$ with the (horizontal) mirror images of our usual notation for $\gamma$ and $\phi$ (respectively); the morphisms $\gamma'$ and $\phi'$ then satisfy the mirror images of the laws in Definition \ref{def:dataserviceaction}.

Furthermore, we can extend the ``filtering'' operation from data services to objects equipped with data service actions.
This enables us to take instances of two entities $M$ and $M'$, each equipped with a data service action by a common data service $D$, and ensure that the corresponding instances of $D$ are equal.
We define this via the morphism $\chi: M \otimes M' \to M \otimes M'$) defined in \eqref{eq:chidef}.

\begin{equation}
\begin{tabular}{cccc}
\begin{tabular}{c}
$\chi :=$
\end{tabular}
&
\begin{tabular}{c}
    \begin{tikzpicture}
        \tikzset{
	% Nodes
	AttrR/.style = {font={\Large$\rhd$}},
	AttrL/.style = {font={\Large$\lhd$}},
	MoveTo/.style = {draw,font={\MoveTo}},
	Pick/.style = {draw,font={\Pick}},
	Place/.style = {draw,font={\Place}},
	Circ/.style = {draw,circle,scale=0.75},
	Rect/.style = {draw,rectangle},
	%Edges
	ob/.style = {out=-90, in=90, looseness=1},
	info/.style = {dashed,out=-90, in=90, looseness=1},
	attrR/.style = {dashed,out=0, in=90, looseness=1},
	attrL/.style = {dashed,out=180, in=90, looseness=1},
	attrRL/.style = {dashed,out=0, in=180, looseness=1},
	attrLR/.style = {dashed,out=180, in=0, looseness=1},
%	attr/.style = {dashed,out=0, in=90, looseness=1},
%	urattr/.style = {dashed,out=-90, in=180, looseness=1},
%	ldattr/.style = {dashed,out=180, in=90, looseness=1},
%	ulattr/.style = {dashed,out=-90, in=0, looseness=1},
%	lrattr/.style = {dashed,out=180, in=0, looseness=1,shorten >=-8pt, shorten <=-8pt},
}

\pgfmathsetmacro{\port}{3mm}

    	% Frame
    	\coordinate (NW) at (0,1); \coordinate (NE) at (3,1); 
    	\coordinate (SW) at (0,0); \coordinate (SE) at (3,0);
    	%\draw[black] (NW) rectangle (SE);
    	
    	%Nodes
    	\node[AttrR] (AttrR1) at (0.9, 0.5) {};
    	\node[AttrL] (AttrL1) at (2.1, 0.5) {};
    	
    	% Autocreate port names
    	\pgfmathsetmacro{\diagramIn}{3}
    	\pgfmathsetmacro{\diagramOut}{3}	
    	\pgfmathsetmacro{\nAttrL}{1}
    	\pgfmathsetmacro{\nAttrR}{1}
    	\pgfmathsetmacro{\nMoveTo}{0}
    	\pgfmathsetmacro{\nPick}{0}
    	\pgfmathsetmacro{\nPlace}{0}
    	\foreach \x [evaluate=\x as \dx using \x/(\diagramIn + 1)] in {1,...,\diagramIn}{
	\coordinate (in\x) at ([yshift=\port]$(NW)!\dx!(NE)$);
}

\foreach \x [evaluate=\x as \dx using \x/(\diagramOut + 1)] in {1,...,\diagramOut}{
	\coordinate (out\x) at ([yshift=-\port]$(SW)!\dx!(SE)$);
}

\ifthenelse	{\nAttrR > 0}
			{\foreach \i in {1,...,\nAttrR}{
				\coordinate (AttrR\i in) at ([xshift=-1.5mm,yshift=1mm]$(AttrR\i)$);
				\coordinate (AttrR\i out) at ([xshift=-1.5mm,yshift=-1mm]$(AttrR\i)$);
				\coordinate (AttrR\i info) at ([xshift=1.5mm,yshift=0mm]$(AttrR\i)$);
			}}{}

\ifthenelse {\nAttrL > 0}
			{\foreach \i in {1,...,\nAttrL}{
				\coordinate (AttrL\i in) at ([xshift=1.5mm,yshift=1mm]$(AttrL\i)$);
				\coordinate (AttrL\i out) at ([xshift=1.5mm,yshift=-1mm]$(AttrL\i)$);
				\coordinate (AttrL\i info) at ([xshift=-1.5mm,yshift=0mm]$(AttrL\i)$);
			}}{}

\ifthenelse {\nMoveTo > 0}
			{\foreach \i in {1,...,\nMoveTo}{
				\coordinate (MoveTo\i inR) at ($(MoveTo\i .north west)!0.33!(MoveTo\i .north east)$);
				\coordinate (MoveTo\i inL) at ($(MoveTo\i .north west)!0.66!(MoveTo\i .north east)$);
				\coordinate (MoveTo\i outR) at ($(MoveTo\i .south west)!0.5!(MoveTo\i .south east)$);
			}}{}

\ifthenelse {\nPick > 0}
			{\foreach \i in {1,...,\nPick}{
					\coordinate (Pick\i inR) at ($(Pick\i .north west)!0.33!(Pick\i .north east)$);
					\coordinate (Pick\i inB) at ($(Pick\i .north west)!0.66!(Pick\i .north east)$);
					\coordinate (Pick\i outR) at ($(Pick\i .south west)!0.5!(Pick\i .south east)$);
			}}{}

\ifthenelse {\nPlace > 0}
			{\foreach \i in {1,...,\nPlace}{
					\coordinate (Place\i inR) at ($(Place\i .north west)!0.5!(Place\i .north east)$);
					\coordinate (Place\i outR) at ($(Place\i .south west)!0.33!(Place\i .south east)$);
					\coordinate (Place\i outB) at ($(Place\i .south west)!0.66!(Place\i .south east)$);
			}}{}
    	
    	% Draw edges
    	\path[every node/.style={font=\scriptsize}]
    	(in1) edge[ob] node[at start, yshift=\port] {$M$} (AttrR1in)
    	(in3) edge[ob] node[at start, yshift=\port] {$M'$} (AttrL1in)
    	(AttrR1) edge[info, out=0, in=180] node[yshift=\port] {$D$} (AttrL1)
    	(AttrR1out) edge[ob] node[at end, yshift=-\port] {$M$} (out1)
    	(AttrL1out) edge[ob] node[at end, yshift=-\port] {$M'$} (out3)
    	;
    \end{tikzpicture}
\end{tabular}
&
\begin{tabular}{c}
:=
\end{tabular}
&
\begin{tabular}{c}
    \begin{tikzpicture}
        \tikzset{
	% Nodes
	AttrR/.style = {font={\Large$\rhd$}},
	AttrL/.style = {font={\Large$\lhd$}},
	MoveTo/.style = {draw,font={\MoveTo}},
	Pick/.style = {draw,font={\Pick}},
	Place/.style = {draw,font={\Place}},
	Circ/.style = {draw,circle,scale=0.75},
	Rect/.style = {draw,rectangle},
	%Edges
	ob/.style = {out=-90, in=90, looseness=1},
	info/.style = {dashed,out=-90, in=90, looseness=1},
	attrR/.style = {dashed,out=0, in=90, looseness=1},
	attrL/.style = {dashed,out=180, in=90, looseness=1},
	attrRL/.style = {dashed,out=0, in=180, looseness=1},
	attrLR/.style = {dashed,out=180, in=0, looseness=1},
%	attr/.style = {dashed,out=0, in=90, looseness=1},
%	urattr/.style = {dashed,out=-90, in=180, looseness=1},
%	ldattr/.style = {dashed,out=180, in=90, looseness=1},
%	ulattr/.style = {dashed,out=-90, in=0, looseness=1},
%	lrattr/.style = {dashed,out=180, in=0, looseness=1,shorten >=-8pt, shorten <=-8pt},
}

\pgfmathsetmacro{\port}{3mm}

    	% Frame
    	\coordinate (NW) at (0,1); \coordinate (NE) at (3,1); 
    	\coordinate (SW) at (0,0); \coordinate (SE) at (3,0);
    	%\draw[black] (NW) rectangle (SE);
    	
    	%Nodes
    	\node[AttrR] (AttrR1) at (0.9, 0.9) {};
    	\node[AttrL] (AttrL1) at (2.1, 0.9) {};
    	\node[Circ] (Circ1) at (1.5, 0.75) {};
    	\node[Circ] (Circ2) at (1.5, 0.1) {};
    	
    	% Autocreate port names
    	\pgfmathsetmacro{\diagramIn}{3}
    	\pgfmathsetmacro{\diagramOut}{3}	
    	\pgfmathsetmacro{\nAttrL}{1}
    	\pgfmathsetmacro{\nAttrR}{1}
    	\pgfmathsetmacro{\nMoveTo}{0}
    	\pgfmathsetmacro{\nPick}{0}
    	\pgfmathsetmacro{\nPlace}{0}
    	\foreach \x [evaluate=\x as \dx using \x/(\diagramIn + 1)] in {1,...,\diagramIn}{
	\coordinate (in\x) at ([yshift=\port]$(NW)!\dx!(NE)$);
}

\foreach \x [evaluate=\x as \dx using \x/(\diagramOut + 1)] in {1,...,\diagramOut}{
	\coordinate (out\x) at ([yshift=-\port]$(SW)!\dx!(SE)$);
}

\ifthenelse	{\nAttrR > 0}
			{\foreach \i in {1,...,\nAttrR}{
				\coordinate (AttrR\i in) at ([xshift=-1.5mm,yshift=1mm]$(AttrR\i)$);
				\coordinate (AttrR\i out) at ([xshift=-1.5mm,yshift=-1mm]$(AttrR\i)$);
				\coordinate (AttrR\i info) at ([xshift=1.5mm,yshift=0mm]$(AttrR\i)$);
			}}{}

\ifthenelse {\nAttrL > 0}
			{\foreach \i in {1,...,\nAttrL}{
				\coordinate (AttrL\i in) at ([xshift=1.5mm,yshift=1mm]$(AttrL\i)$);
				\coordinate (AttrL\i out) at ([xshift=1.5mm,yshift=-1mm]$(AttrL\i)$);
				\coordinate (AttrL\i info) at ([xshift=-1.5mm,yshift=0mm]$(AttrL\i)$);
			}}{}

\ifthenelse {\nMoveTo > 0}
			{\foreach \i in {1,...,\nMoveTo}{
				\coordinate (MoveTo\i inR) at ($(MoveTo\i .north west)!0.33!(MoveTo\i .north east)$);
				\coordinate (MoveTo\i inL) at ($(MoveTo\i .north west)!0.66!(MoveTo\i .north east)$);
				\coordinate (MoveTo\i outR) at ($(MoveTo\i .south west)!0.5!(MoveTo\i .south east)$);
			}}{}

\ifthenelse {\nPick > 0}
			{\foreach \i in {1,...,\nPick}{
					\coordinate (Pick\i inR) at ($(Pick\i .north west)!0.33!(Pick\i .north east)$);
					\coordinate (Pick\i inB) at ($(Pick\i .north west)!0.66!(Pick\i .north east)$);
					\coordinate (Pick\i outR) at ($(Pick\i .south west)!0.5!(Pick\i .south east)$);
			}}{}

\ifthenelse {\nPlace > 0}
			{\foreach \i in {1,...,\nPlace}{
					\coordinate (Place\i inR) at ($(Place\i .north west)!0.5!(Place\i .north east)$);
					\coordinate (Place\i outR) at ($(Place\i .south west)!0.33!(Place\i .south east)$);
					\coordinate (Place\i outB) at ($(Place\i .south west)!0.66!(Place\i .south east)$);
			}}{}
    	
    	% Draw edges
    	\path[every node/.style={font=\scriptsize}]
    	(in1) edge[ob] node[at start, yshift=\port] {$M$} (AttrR1in)
    	(in3) edge[ob] node[at start, yshift=\port] {$M'$} (AttrL1in)
    	(AttrR1) edge[info, out=0, in=180] node[yshift=\port] {$D$} (Circ1)
    	(AttrL1) edge[info, out=180, in=0] node[yshift=\port] {$D$} (Circ1)
    	(AttrR1out) edge[ob] node[at end, yshift=-\port] {$M$} (out1)
    	(AttrL1out) edge[ob] node[at end, yshift=-\port] {$M'$} (out3)
    	(Circ1) edge[info] node[left] {$D$} (Circ2)
    	;
    \end{tikzpicture}
\end{tabular}
\end{tabular}\end{equation}

Our depiction of $\chi$ is \emph{a priori} ambiguous, as it could easily be interpreted as referring to other morphisms, for example:
\begin{equation} \label{eq:chidef}
\begin{tabular}{ccc}
\begin{tabular}{c}
    \begin{tikzpicture}
        \tikzset{
	% Nodes
	AttrR/.style = {font={\Large$\rhd$}},
	AttrL/.style = {font={\Large$\lhd$}},
	MoveTo/.style = {draw,font={\MoveTo}},
	Pick/.style = {draw,font={\Pick}},
	Place/.style = {draw,font={\Place}},
	Circ/.style = {draw,circle,scale=0.75},
	Rect/.style = {draw,rectangle},
	%Edges
	ob/.style = {out=-90, in=90, looseness=1},
	info/.style = {dashed,out=-90, in=90, looseness=1},
	attrR/.style = {dashed,out=0, in=90, looseness=1},
	attrL/.style = {dashed,out=180, in=90, looseness=1},
	attrRL/.style = {dashed,out=0, in=180, looseness=1},
	attrLR/.style = {dashed,out=180, in=0, looseness=1},
%	attr/.style = {dashed,out=0, in=90, looseness=1},
%	urattr/.style = {dashed,out=-90, in=180, looseness=1},
%	ldattr/.style = {dashed,out=180, in=90, looseness=1},
%	ulattr/.style = {dashed,out=-90, in=0, looseness=1},
%	lrattr/.style = {dashed,out=180, in=0, looseness=1,shorten >=-8pt, shorten <=-8pt},
}

\pgfmathsetmacro{\port}{3mm}

    	% Frame
    	\coordinate (NW) at (0,1); \coordinate (NE) at (3,1); 
    	\coordinate (SW) at (0,0); \coordinate (SE) at (3,0);
    	%\draw[black] (NW) rectangle (SE);
    	
    	%Nodes
    	\node[AttrR] (AttrR1) at (0.9, 0.5) {};
    	\node[AttrL] (AttrL1) at (2.1, 0.5) {};
    	
    	% Autocreate port names
    	\pgfmathsetmacro{\diagramIn}{3}
    	\pgfmathsetmacro{\diagramOut}{3}	
    	\pgfmathsetmacro{\nAttrL}{1}
    	\pgfmathsetmacro{\nAttrR}{1}
    	\pgfmathsetmacro{\nMoveTo}{0}
    	\pgfmathsetmacro{\nPick}{0}
    	\pgfmathsetmacro{\nPlace}{0}
    	\foreach \x [evaluate=\x as \dx using \x/(\diagramIn + 1)] in {1,...,\diagramIn}{
	\coordinate (in\x) at ([yshift=\port]$(NW)!\dx!(NE)$);
}

\foreach \x [evaluate=\x as \dx using \x/(\diagramOut + 1)] in {1,...,\diagramOut}{
	\coordinate (out\x) at ([yshift=-\port]$(SW)!\dx!(SE)$);
}

\ifthenelse	{\nAttrR > 0}
			{\foreach \i in {1,...,\nAttrR}{
				\coordinate (AttrR\i in) at ([xshift=-1.5mm,yshift=1mm]$(AttrR\i)$);
				\coordinate (AttrR\i out) at ([xshift=-1.5mm,yshift=-1mm]$(AttrR\i)$);
				\coordinate (AttrR\i info) at ([xshift=1.5mm,yshift=0mm]$(AttrR\i)$);
			}}{}

\ifthenelse {\nAttrL > 0}
			{\foreach \i in {1,...,\nAttrL}{
				\coordinate (AttrL\i in) at ([xshift=1.5mm,yshift=1mm]$(AttrL\i)$);
				\coordinate (AttrL\i out) at ([xshift=1.5mm,yshift=-1mm]$(AttrL\i)$);
				\coordinate (AttrL\i info) at ([xshift=-1.5mm,yshift=0mm]$(AttrL\i)$);
			}}{}

\ifthenelse {\nMoveTo > 0}
			{\foreach \i in {1,...,\nMoveTo}{
				\coordinate (MoveTo\i inR) at ($(MoveTo\i .north west)!0.33!(MoveTo\i .north east)$);
				\coordinate (MoveTo\i inL) at ($(MoveTo\i .north west)!0.66!(MoveTo\i .north east)$);
				\coordinate (MoveTo\i outR) at ($(MoveTo\i .south west)!0.5!(MoveTo\i .south east)$);
			}}{}

\ifthenelse {\nPick > 0}
			{\foreach \i in {1,...,\nPick}{
					\coordinate (Pick\i inR) at ($(Pick\i .north west)!0.33!(Pick\i .north east)$);
					\coordinate (Pick\i inB) at ($(Pick\i .north west)!0.66!(Pick\i .north east)$);
					\coordinate (Pick\i outR) at ($(Pick\i .south west)!0.5!(Pick\i .south east)$);
			}}{}

\ifthenelse {\nPlace > 0}
			{\foreach \i in {1,...,\nPlace}{
					\coordinate (Place\i inR) at ($(Place\i .north west)!0.5!(Place\i .north east)$);
					\coordinate (Place\i outR) at ($(Place\i .south west)!0.33!(Place\i .south east)$);
					\coordinate (Place\i outB) at ($(Place\i .south west)!0.66!(Place\i .south east)$);
			}}{}
    	
    	% Draw edges
    	\path[every node/.style={font=\scriptsize}]
    	(in1) edge[ob] node[at start, yshift=\port] {$M$} (AttrR1in)
    	(in3) edge[ob] node[at start, yshift=\port] {$M'$} (AttrL1in)
    	(AttrR1) edge[info, out=0, in=180] node[yshift=\port] {$D$} (AttrL1)
    	(AttrR1out) edge[ob] node[at end, yshift=-\port] {$M$} (out1)
    	(AttrL1out) edge[ob] node[at end, yshift=-\port] {$M'$} (out3)
    	;
    \end{tikzpicture}
\end{tabular}
&
\begin{tabular}{c}
$\overset{?}{=}$
\end{tabular}
&
\begin{tabular}{c}
    \begin{tikzpicture}
        \tikzset{
	% Nodes
	AttrR/.style = {font={\Large$\rhd$}},
	AttrL/.style = {font={\Large$\lhd$}},
	MoveTo/.style = {draw,font={\MoveTo}},
	Pick/.style = {draw,font={\Pick}},
	Place/.style = {draw,font={\Place}},
	Circ/.style = {draw,circle,scale=0.75},
	Rect/.style = {draw,rectangle},
	%Edges
	ob/.style = {out=-90, in=90, looseness=1},
	info/.style = {dashed,out=-90, in=90, looseness=1},
	attrR/.style = {dashed,out=0, in=90, looseness=1},
	attrL/.style = {dashed,out=180, in=90, looseness=1},
	attrRL/.style = {dashed,out=0, in=180, looseness=1},
	attrLR/.style = {dashed,out=180, in=0, looseness=1},
%	attr/.style = {dashed,out=0, in=90, looseness=1},
%	urattr/.style = {dashed,out=-90, in=180, looseness=1},
%	ldattr/.style = {dashed,out=180, in=90, looseness=1},
%	ulattr/.style = {dashed,out=-90, in=0, looseness=1},
%	lrattr/.style = {dashed,out=180, in=0, looseness=1,shorten >=-8pt, shorten <=-8pt},
}

\pgfmathsetmacro{\port}{3mm}

    	% Frame
    	\coordinate (NW) at (0,1); \coordinate (NE) at (3,1); 
    	\coordinate (SW) at (0,0); \coordinate (SE) at (3,0);
    	%\draw[black] (NW) rectangle (SE);
    	
    	%Nodes
    	\node[AttrR] (AttrR1) at (0.9, 0.25) {};
    	\node[AttrL] (AttrL1) at (2.1, 0.75) {};
    	
    	% Autocreate port names
    	\pgfmathsetmacro{\diagramIn}{3}
    	\pgfmathsetmacro{\diagramOut}{3}	
    	\pgfmathsetmacro{\nAttrL}{1}
    	\pgfmathsetmacro{\nAttrR}{1}
    	\pgfmathsetmacro{\nMoveTo}{0}
    	\pgfmathsetmacro{\nPick}{0}
    	\pgfmathsetmacro{\nPlace}{0}
    	\foreach \x [evaluate=\x as \dx using \x/(\diagramIn + 1)] in {1,...,\diagramIn}{
	\coordinate (in\x) at ([yshift=\port]$(NW)!\dx!(NE)$);
}

\foreach \x [evaluate=\x as \dx using \x/(\diagramOut + 1)] in {1,...,\diagramOut}{
	\coordinate (out\x) at ([yshift=-\port]$(SW)!\dx!(SE)$);
}

\ifthenelse	{\nAttrR > 0}
			{\foreach \i in {1,...,\nAttrR}{
				\coordinate (AttrR\i in) at ([xshift=-1.5mm,yshift=1mm]$(AttrR\i)$);
				\coordinate (AttrR\i out) at ([xshift=-1.5mm,yshift=-1mm]$(AttrR\i)$);
				\coordinate (AttrR\i info) at ([xshift=1.5mm,yshift=0mm]$(AttrR\i)$);
			}}{}

\ifthenelse {\nAttrL > 0}
			{\foreach \i in {1,...,\nAttrL}{
				\coordinate (AttrL\i in) at ([xshift=1.5mm,yshift=1mm]$(AttrL\i)$);
				\coordinate (AttrL\i out) at ([xshift=1.5mm,yshift=-1mm]$(AttrL\i)$);
				\coordinate (AttrL\i info) at ([xshift=-1.5mm,yshift=0mm]$(AttrL\i)$);
			}}{}

\ifthenelse {\nMoveTo > 0}
			{\foreach \i in {1,...,\nMoveTo}{
				\coordinate (MoveTo\i inR) at ($(MoveTo\i .north west)!0.33!(MoveTo\i .north east)$);
				\coordinate (MoveTo\i inL) at ($(MoveTo\i .north west)!0.66!(MoveTo\i .north east)$);
				\coordinate (MoveTo\i outR) at ($(MoveTo\i .south west)!0.5!(MoveTo\i .south east)$);
			}}{}

\ifthenelse {\nPick > 0}
			{\foreach \i in {1,...,\nPick}{
					\coordinate (Pick\i inR) at ($(Pick\i .north west)!0.33!(Pick\i .north east)$);
					\coordinate (Pick\i inB) at ($(Pick\i .north west)!0.66!(Pick\i .north east)$);
					\coordinate (Pick\i outR) at ($(Pick\i .south west)!0.5!(Pick\i .south east)$);
			}}{}

\ifthenelse {\nPlace > 0}
			{\foreach \i in {1,...,\nPlace}{
					\coordinate (Place\i inR) at ($(Place\i .north west)!0.5!(Place\i .north east)$);
					\coordinate (Place\i outR) at ($(Place\i .south west)!0.33!(Place\i .south east)$);
					\coordinate (Place\i outB) at ($(Place\i .south west)!0.66!(Place\i .south east)$);
			}}{}
    	
    	% Draw edges
    	\path[every node/.style={font=\scriptsize}]
    	(in1) edge[ob] node[at start, yshift=\port] {$M$} (AttrR1in)
    	(in3) edge[ob] node[at start, yshift=\port] {$M'$} (AttrL1in)
    	(AttrR1) edge[info, out=0, in=180] node[yshift=\port] {$D$} (AttrL1)
    	(AttrR1out) edge[ob] node[at end, yshift=-\port] {$M$} (out1)
    	(AttrL1out) edge[ob] node[at end, yshift=-\port] {$M'$} (out3)
    	;
    \end{tikzpicture}
\end{tabular}
\end{tabular}\end{equation}

However, there is no real ambiguity here, as in fact this equation (and other obvious equations arising from our depiction) do hold, e.g.\ by Proposition \ref{prop:comonoid}.
Proposition \ref{prop:comonoid} can also be used to show that the diagrams in Figure \ref{fig:proofs} (and other well-formed diagrams with the action morphisms displayed horizontally) represent well-defined morphisms.

From this setup, the definition of a category with attributes (or equivalently an attribute structure on a SMC) is straightforward.

\begin{defn} \label{def:attributestructure}
An \emph{attribute structure} $(\AA, E, V, \gamma)$ on a SMC $\CC$ consists of:
\begin{itemize}
    \item A category $\AA$, called the \emph{category of attributes};
    \item A functor $E: \AA \to \CC$;
    \item A functor $V: \AA \to \Data(\CC)$;
    \item A natural transformation\footnote{Recall $U: \Data(\CC) \to \CC$ is the forgetful functor.} $\gamma: E \to E \otimes (U \circ V)$ 
\end{itemize}
such that:
\begin{itemize}
    \item Each $\gamma_A: E(A) \to E(A) \otimes V(A)$ gives a action of the comonoid $(V(A), \delta_{V(A)}, \epsilon_{V(A)})$ on $E(A)$.
\end{itemize}
When we have an attribute structure on $\CC$ in mind, we refer to $\CC$ as a \emph{category with attributes}.
\end{defn}

An attribute structure on $\CC$ be viewed as a ``collection of distinguished attributes in $\CC$.''
The naturality condition on $\gamma$ ensures that, given $f: A \to A'$ in $\AA$, the map $V(f)$ allows us to ``predict'' the effects of $E(f)$ on the relevant informational objects.
Note that we require no naturality condition on the $\phi_A$ maps associated to each $\gamma_A$ through Proposition \ref{prop:comonoid}.
This is because we cannot anticipate that the ``processing'' performed by $V(f)$ will produce meaningful results when applied to an arbitrary input.

\subsection{$\Poset$-Enriched Case} \label{sec:posetcase}

When modeling attributes, we are often interested in SMCs which are enriched in $\Poset$.
This enriched structure can be used to model the notion that a morphism or process is partially defined, i.e.\ that it could fail to execute on certain inputs.
Heuristically, morphisms $f$ and $g$ satisfy $f \leq g$ iff the results of $f$ are the same as those of $g$ when both exist, and whenever $f$ successfully executes given a certain input, so too does $g$.
We spell out the details of this enrichment for ease of reference.

% WE PROBABLY DON'T NEED TO SPELL THINGS OUT.
% BUT WE SHOULD AT LEAST MAKE CLEAR HOW THE POSET STRUCTURE INTERACTS WITH TENSOR PRODUCTS

\begin{defn}
A $\Poset$-enrichment of a category $\CC$ is an assignment of a partial order (uniformly denoted $\leq$) to each hom-set $\CC(A, B)$ such that:
\begin{itemize}
    \item If $f \leq f'$ and $g \leq g'$, and $g \circ f$ and $g' \circ f'$ exist, then $g \circ f \leq g' \circ f'$.
\end{itemize}
If $(\CC, \otimes, I)$ is a SMC, we require furthermore:
\begin{itemize}
\item If $f \leq f'$ and $g \leq g'$, then $f \otimes g \leq f' \otimes g'$.
\end{itemize}
A $\Poset$-enriched functor between $\Poset$-enriched categories $\CC, \DD$ is a functor $F: \CC \to \DD$ such that the associated functions $\CC(A, B) \to \DD(FA, FB)$ are all monotone increasing.
\end{defn}

From this viewpoint, a certain condition on the data services in a $\Poset$-enriched SMC seems desirable.
Namely, filtering two values for equality and then returning both values should yield the same results as simply returning both values whenever the former operation is defined.
This is best formalized in the following axiom, which appears in \cite{bonchi2017functorial}.

\begin{defn}
A data service $D$ in a $\Poset$-enriched SMC $\CC$ is said to be \emph{well-behaved with respect to the enrichment} (for brevity, ``\emph{well-behaved}'') if it satisfies:
\begin{equation}
\end{equation} 
\end{defn}

One is led to wonder whether or not an analogous property to that satisfied by well-behaved data services holds for data service actions by such data services.
This is indeed the case, as shown in the following proposition.

\begin{prop} \label{prop:posetineq1}
Let $\CC$ be a $\Poset$-enriched SMC, $D$ a well-behaved data service in $\CC$, and $M$ an object of $\CC$ together with a data service action by $D$. Then:
\begin{equation}%
\end{equation}
\end{prop}

We can also prove a useful proposition about the morphism $\chi$ described above.

\begin{prop} \label{prop:posetineq2}
Let $\CC$ be a $\Poset$-enriched SMC, $D$ a well-behaved data service in $\CC$, and $M$ an object of $\CC$ together with a data service action by $D$. Then:
\begin{equation}%
\end{tabular}\end{equation}
\end{prop}

As indicated above, well-behaved data services enable us to create a useful variation on the category $\Data(\CC)$ defined previously.

\begin{defn} \label{def:dataprime}
Let $(\CC, \otimes, I, \leq)$ be a $\Poset$-enriched SMC. Define $\Data'(\CC)$, the \emph{category of well-behaved data services} in $\CC$, to be the category where:
\begin{itemize}
    \item Objects are well-behaved data services in $\CC$, and
    \item Morphisms $D \to D'$ are \emph{lax data service homomorphisms} from $D$ to $D'$; that is, morphisms $f: D \to D'$ satisfying the axioms (see \cite{bonchi2017functorial}):
\end{itemize}
\begin{equation} \label{eq:weakhomcomult} 
\begin{tabular}{ccc}
\begin{tabular}{c}
    \begin{tikzpicture}
        \tikzset{
	% Nodes
	AttrR/.style = {font={\Large$\rhd$}},
	AttrL/.style = {font={\Large$\lhd$}},
	MoveTo/.style = {draw,font={\MoveTo}},
	Pick/.style = {draw,font={\Pick}},
	Place/.style = {draw,font={\Place}},
	Circ/.style = {draw,circle,scale=0.75},
	Rect/.style = {draw,rectangle},
	%Edges
	ob/.style = {out=-90, in=90, looseness=1},
	info/.style = {dashed,out=-90, in=90, looseness=1},
	attrR/.style = {dashed,out=0, in=90, looseness=1},
	attrL/.style = {dashed,out=180, in=90, looseness=1},
	attrRL/.style = {dashed,out=0, in=180, looseness=1},
	attrLR/.style = {dashed,out=180, in=0, looseness=1},
%	attr/.style = {dashed,out=0, in=90, looseness=1},
%	urattr/.style = {dashed,out=-90, in=180, looseness=1},
%	ldattr/.style = {dashed,out=180, in=90, looseness=1},
%	ulattr/.style = {dashed,out=-90, in=0, looseness=1},
%	lrattr/.style = {dashed,out=180, in=0, looseness=1,shorten >=-8pt, shorten <=-8pt},
}

\pgfmathsetmacro{\port}{3mm}

    	% Frame
    	\coordinate (NW) at (0,1); \coordinate (NE) at (3,1); 
    	\coordinate (SW) at (0,0); \coordinate (SE) at (3,0);
    	%\draw[black] (NW) rectangle (SE);
    	
    	%Nodes
    	\node[Rect] (Rect1) at (1.5, 0.9) {$f$};
    	\node[Circ] (Circ1) at (1.5, 0.1) {};
    	
    	% Autocreate port names
    	\pgfmathsetmacro{\diagramIn}{1}
    	\pgfmathsetmacro{\diagramOut}{2}	
    	\pgfmathsetmacro{\nAttrL}{0}
    	\pgfmathsetmacro{\nAttrR}{0}
    	\pgfmathsetmacro{\nMoveTo}{0}
    	\pgfmathsetmacro{\nPick}{0}
    	\pgfmathsetmacro{\nPlace}{0}
    	\foreach \x [evaluate=\x as \dx using \x/(\diagramIn + 1)] in {1,...,\diagramIn}{
	\coordinate (in\x) at ([yshift=\port]$(NW)!\dx!(NE)$);
}

\foreach \x [evaluate=\x as \dx using \x/(\diagramOut + 1)] in {1,...,\diagramOut}{
	\coordinate (out\x) at ([yshift=-\port]$(SW)!\dx!(SE)$);
}

\ifthenelse	{\nAttrR > 0}
			{\foreach \i in {1,...,\nAttrR}{
				\coordinate (AttrR\i in) at ([xshift=-1.5mm,yshift=1mm]$(AttrR\i)$);
				\coordinate (AttrR\i out) at ([xshift=-1.5mm,yshift=-1mm]$(AttrR\i)$);
				\coordinate (AttrR\i info) at ([xshift=1.5mm,yshift=0mm]$(AttrR\i)$);
			}}{}

\ifthenelse {\nAttrL > 0}
			{\foreach \i in {1,...,\nAttrL}{
				\coordinate (AttrL\i in) at ([xshift=1.5mm,yshift=1mm]$(AttrL\i)$);
				\coordinate (AttrL\i out) at ([xshift=1.5mm,yshift=-1mm]$(AttrL\i)$);
				\coordinate (AttrL\i info) at ([xshift=-1.5mm,yshift=0mm]$(AttrL\i)$);
			}}{}

\ifthenelse {\nMoveTo > 0}
			{\foreach \i in {1,...,\nMoveTo}{
				\coordinate (MoveTo\i inR) at ($(MoveTo\i .north west)!0.33!(MoveTo\i .north east)$);
				\coordinate (MoveTo\i inL) at ($(MoveTo\i .north west)!0.66!(MoveTo\i .north east)$);
				\coordinate (MoveTo\i outR) at ($(MoveTo\i .south west)!0.5!(MoveTo\i .south east)$);
			}}{}

\ifthenelse {\nPick > 0}
			{\foreach \i in {1,...,\nPick}{
					\coordinate (Pick\i inR) at ($(Pick\i .north west)!0.33!(Pick\i .north east)$);
					\coordinate (Pick\i inB) at ($(Pick\i .north west)!0.66!(Pick\i .north east)$);
					\coordinate (Pick\i outR) at ($(Pick\i .south west)!0.5!(Pick\i .south east)$);
			}}{}

\ifthenelse {\nPlace > 0}
			{\foreach \i in {1,...,\nPlace}{
					\coordinate (Place\i inR) at ($(Place\i .north west)!0.5!(Place\i .north east)$);
					\coordinate (Place\i outR) at ($(Place\i .south west)!0.33!(Place\i .south east)$);
					\coordinate (Place\i outB) at ($(Place\i .south west)!0.66!(Place\i .south east)$);
			}}{}
    	
    	% Draw edges
    	\path[every node/.style={font=\scriptsize}]
    	(in1) edge[info] node[at start, yshift=\port] {$D$} (Rect1)
    	(Rect1) edge[info] node[left] {$D'$} (Circ1)
    	(Circ1) edge[info, out=180] node[at end, yshift=-\port] {$D'$} (out1)
    	(Circ1) edge[info, out=0] node[at end, yshift=-\port] {$D'$} (out2)
    	;
    \end{tikzpicture}
\end{tabular}
&
\begin{tabular}{c}
$\leq$
\end{tabular}
&
\begin{tabular}{c}
    \begin{tikzpicture}
        \tikzset{
	% Nodes
	AttrR/.style = {font={\Large$\rhd$}},
	AttrL/.style = {font={\Large$\lhd$}},
	MoveTo/.style = {draw,font={\MoveTo}},
	Pick/.style = {draw,font={\Pick}},
	Place/.style = {draw,font={\Place}},
	Circ/.style = {draw,circle,scale=0.75},
	Rect/.style = {draw,rectangle},
	%Edges
	ob/.style = {out=-90, in=90, looseness=1},
	info/.style = {dashed,out=-90, in=90, looseness=1},
	attrR/.style = {dashed,out=0, in=90, looseness=1},
	attrL/.style = {dashed,out=180, in=90, looseness=1},
	attrRL/.style = {dashed,out=0, in=180, looseness=1},
	attrLR/.style = {dashed,out=180, in=0, looseness=1},
%	attr/.style = {dashed,out=0, in=90, looseness=1},
%	urattr/.style = {dashed,out=-90, in=180, looseness=1},
%	ldattr/.style = {dashed,out=180, in=90, looseness=1},
%	ulattr/.style = {dashed,out=-90, in=0, looseness=1},
%	lrattr/.style = {dashed,out=180, in=0, looseness=1,shorten >=-8pt, shorten <=-8pt},
}

\pgfmathsetmacro{\port}{3mm}

    	% Frame
    	\coordinate (NW) at (0,1); \coordinate (NE) at (3,1); 
    	\coordinate (SW) at (0,0); \coordinate (SE) at (3,0);
    	%\draw[black] (NW) rectangle (SE);
    	
    	%Nodes
    	\node[Rect] (Rect1) at (1.0, 0.1) {$f$};
    	\node[Rect] (Rect2) at (2.0, 0.1) {$f$};
    	\node[Circ] (Circ1) at (1.5, 0.9) {};
    	
    	% Autocreate port names
    	\pgfmathsetmacro{\diagramIn}{1}
    	\pgfmathsetmacro{\diagramOut}{2}	
    	\pgfmathsetmacro{\nAttrL}{0}
    	\pgfmathsetmacro{\nAttrR}{0}
    	\pgfmathsetmacro{\nMoveTo}{0}
    	\pgfmathsetmacro{\nPick}{0}
    	\pgfmathsetmacro{\nPlace}{0}
    	\foreach \x [evaluate=\x as \dx using \x/(\diagramIn + 1)] in {1,...,\diagramIn}{
	\coordinate (in\x) at ([yshift=\port]$(NW)!\dx!(NE)$);
}

\foreach \x [evaluate=\x as \dx using \x/(\diagramOut + 1)] in {1,...,\diagramOut}{
	\coordinate (out\x) at ([yshift=-\port]$(SW)!\dx!(SE)$);
}

\ifthenelse	{\nAttrR > 0}
			{\foreach \i in {1,...,\nAttrR}{
				\coordinate (AttrR\i in) at ([xshift=-1.5mm,yshift=1mm]$(AttrR\i)$);
				\coordinate (AttrR\i out) at ([xshift=-1.5mm,yshift=-1mm]$(AttrR\i)$);
				\coordinate (AttrR\i info) at ([xshift=1.5mm,yshift=0mm]$(AttrR\i)$);
			}}{}

\ifthenelse {\nAttrL > 0}
			{\foreach \i in {1,...,\nAttrL}{
				\coordinate (AttrL\i in) at ([xshift=1.5mm,yshift=1mm]$(AttrL\i)$);
				\coordinate (AttrL\i out) at ([xshift=1.5mm,yshift=-1mm]$(AttrL\i)$);
				\coordinate (AttrL\i info) at ([xshift=-1.5mm,yshift=0mm]$(AttrL\i)$);
			}}{}

\ifthenelse {\nMoveTo > 0}
			{\foreach \i in {1,...,\nMoveTo}{
				\coordinate (MoveTo\i inR) at ($(MoveTo\i .north west)!0.33!(MoveTo\i .north east)$);
				\coordinate (MoveTo\i inL) at ($(MoveTo\i .north west)!0.66!(MoveTo\i .north east)$);
				\coordinate (MoveTo\i outR) at ($(MoveTo\i .south west)!0.5!(MoveTo\i .south east)$);
			}}{}

\ifthenelse {\nPick > 0}
			{\foreach \i in {1,...,\nPick}{
					\coordinate (Pick\i inR) at ($(Pick\i .north west)!0.33!(Pick\i .north east)$);
					\coordinate (Pick\i inB) at ($(Pick\i .north west)!0.66!(Pick\i .north east)$);
					\coordinate (Pick\i outR) at ($(Pick\i .south west)!0.5!(Pick\i .south east)$);
			}}{}

\ifthenelse {\nPlace > 0}
			{\foreach \i in {1,...,\nPlace}{
					\coordinate (Place\i inR) at ($(Place\i .north west)!0.5!(Place\i .north east)$);
					\coordinate (Place\i outR) at ($(Place\i .south west)!0.33!(Place\i .south east)$);
					\coordinate (Place\i outB) at ($(Place\i .south west)!0.66!(Place\i .south east)$);
			}}{}
    	
    	% Draw edges
    	\path[every node/.style={font=\scriptsize}]
    	(in1) edge[info] node[at start, yshift=\port] {$D$} (Circ1)
    	(Circ1) edge[info, out=180] node[left] {$D$} (Rect1)
    	(Circ1) edge[info, out=0] node[right] {$D$} (Rect2)
    	(Rect1) edge[info] node[at end, yshift=-\port] {$D'$} (out1)
    	(Rect2) edge[info] node[at end, yshift=-\port] {$D'$} (out2)
    	;
    \end{tikzpicture}
\end{tabular}
\end{tabular}
\end{equation}

\begin{equation} \label{eq:weakhomcounit}
\begin{tabular}{ccc}
\begin{tabular}{c}
    \begin{tikzpicture}
        \tikzset{
	% Nodes
	AttrR/.style = {font={\Large$\rhd$}},
	AttrL/.style = {font={\Large$\lhd$}},
	MoveTo/.style = {draw,font={\MoveTo}},
	Pick/.style = {draw,font={\Pick}},
	Place/.style = {draw,font={\Place}},
	Circ/.style = {draw,circle,scale=0.75},
	Rect/.style = {draw,rectangle},
	%Edges
	ob/.style = {out=-90, in=90, looseness=1},
	info/.style = {dashed,out=-90, in=90, looseness=1},
	attrR/.style = {dashed,out=0, in=90, looseness=1},
	attrL/.style = {dashed,out=180, in=90, looseness=1},
	attrRL/.style = {dashed,out=0, in=180, looseness=1},
	attrLR/.style = {dashed,out=180, in=0, looseness=1},
%	attr/.style = {dashed,out=0, in=90, looseness=1},
%	urattr/.style = {dashed,out=-90, in=180, looseness=1},
%	ldattr/.style = {dashed,out=180, in=90, looseness=1},
%	ulattr/.style = {dashed,out=-90, in=0, looseness=1},
%	lrattr/.style = {dashed,out=180, in=0, looseness=1,shorten >=-8pt, shorten <=-8pt},
}

\pgfmathsetmacro{\port}{3mm}

    	% Frame
    	\coordinate (NW) at (0,1); \coordinate (NE) at (3,1); 
    	\coordinate (SW) at (0,0); \coordinate (SE) at (3,0);
    	%\draw[black] (NW) rectangle (SE);
    	
    	%Nodes
    	\node[Rect] (Rect1) at (1.5, 0.9) {$f$};
    	\node[Circ] (Circ1) at (1.5, 0.1) {};
    	
    	% Autocreate port names
    	\pgfmathsetmacro{\diagramIn}{1}
    	\pgfmathsetmacro{\diagramOut}{0}	
    	\pgfmathsetmacro{\nAttrL}{0}
    	\pgfmathsetmacro{\nAttrR}{0}
    	\pgfmathsetmacro{\nMoveTo}{0}
    	\pgfmathsetmacro{\nPick}{0}
    	\pgfmathsetmacro{\nPlace}{0}

    	% Draw edges
    	\path[every node/.style={font=\scriptsize}]
    	(in1) edge[info] node[at start, yshift=\port] {$D$} (Rect1)
    	(Rect1) edge[info] node[left] {$D'$} (Circ1)
    	;
    \end{tikzpicture}
\end{tabular}
&
\begin{tabular}{c}
$\leq$
\end{tabular}
&
\begin{tabular}{c}
    \begin{tikzpicture}
        
    	% Frame
    	\coordinate (NW) at (0,1); \coordinate (NE) at (3,1); 
    	\coordinate (SW) at (0,0); \coordinate (SE) at (3,0);
    	%\draw[black] (NW) rectangle (SE);
    	
    	%Nodes
    	\node[Circ] (Circ1) at (1.5, 0.1) {};
    	
    	% Autocreate port names
    	\pgfmathsetmacro{\diagramIn}{1}
    	\pgfmathsetmacro{\diagramOut}{2}	
    	\pgfmathsetmacro{\nAttrL}{0}
    	\pgfmathsetmacro{\nAttrR}{0}
    	\pgfmathsetmacro{\nMoveTo}{0}
    	\pgfmathsetmacro{\nPick}{0}
    	\pgfmathsetmacro{\nPlace}{0}

    	% Draw edges
    	\path[every node/.style={font=\scriptsize}]
    	(in1) edge[info] node[at start, yshift=\port] {$D$} (Circ1)
    	;
    \end{tikzpicture}
\end{tabular}
\end{tabular}
\end{equation}

\begin{equation} \label{eq:weakhommult}
\begin{tabular}{ccc}
\begin{tabular}{c}
    \begin{tikzpicture}
        
    	% Frame
    	\coordinate (NW) at (0,1); \coordinate (NE) at (3,1); 
    	\coordinate (SW) at (0,0); \coordinate (SE) at (3,0);
    	%\draw[black] (NW) rectangle (SE);
    	
    	%Nodes
    	\node[Circ] (Circ1) at (1.5, 0.9) {};
    	\node[Rect] (Rect1) at (1.5, 0.1) {$f$};
    	
    	% Autocreate port names
    	\pgfmathsetmacro{\diagramIn}{2}
    	\pgfmathsetmacro{\diagramOut}{1}	
    	\pgfmathsetmacro{\nAttrL}{0}
    	\pgfmathsetmacro{\nAttrR}{0}
    	\pgfmathsetmacro{\nMoveTo}{0}
    	\pgfmathsetmacro{\nPick}{0}
    	\pgfmathsetmacro{\nPlace}{0}

    	% Draw edges
    	\path[every node/.style={font=\scriptsize}]
    	(in1) edge[info, in=180] node[at start, yshift=\port] {$D$} (Circ1)
    	(in2) edge[info, in=0] node[at start, yshift=\port] {$D$} (Circ1)
    	(Circ1) edge[info] node[left] {$D$} (Rect1)
    	(Rect1) edge[info] node[at end, yshift=-\port] {$D'$} (out1)
    	;
    \end{tikzpicture}
\end{tabular}
&
\begin{tabular}{c}
$\leq$
\end{tabular}
&
\begin{tabular}{c}
    \begin{tikzpicture}
        
    	% Frame
    	\coordinate (NW) at (0,1); \coordinate (NE) at (3,1); 
    	\coordinate (SW) at (0,0); \coordinate (SE) at (3,0);
    	%\draw[black] (NW) rectangle (SE);
    	
    	%Nodes
    	\node[Circ] (Circ1) at (1.5, 0.1) {};
    	\node[Rect] (Rect1) at (1.0, 0.9) {$f$};
    	\node[Rect] (Rect2) at (2.0, 0.9) {$f$};
    	
    	% Autocreate port names
    	\pgfmathsetmacro{\diagramIn}{2}
    	\pgfmathsetmacro{\diagramOut}{1}	
    	\pgfmathsetmacro{\nAttrL}{0}
    	\pgfmathsetmacro{\nAttrR}{0}
    	\pgfmathsetmacro{\nMoveTo}{0}
    	\pgfmathsetmacro{\nPick}{0}
    	\pgfmathsetmacro{\nPlace}{0}

    	% Draw edges
    	\path[every node/.style={font=\scriptsize}]
    	(in1) edge[info] node[at start, yshift=\port] {$D$} (Rect1)
    	(in2) edge[info] node[at start, yshift=\port] {$D$} (Rect2)
    	(Rect1) edge[info, in=180] node[left] {$D'$} (Circ1)
    	(Rect2) edge[info, in=0] node[right] {$D'$} (Circ1)
    	(Circ1) edge[info] node[at end, yshift=-\port] {$D'$} (out1)
    	;
    \end{tikzpicture}
\end{tabular}
\end{tabular}\end{equation}

A morphism $f$ in $\Data'(\CC)$ is \emph{deterministic} if it is a (strong) comonoid homomorphism, i.e.\ if both \eqref{eq:weakhomcomult} and \eqref{eq:weakhomcounit} are equalities.
\end{defn}

The inequalities in Definition \ref{def:dataprime} are described in \cite{bonchi2017functorial} as (some of) those which hold between relations and ($\Set$-based) monoids.
That is, each object of the categories $\Rel$ or $\PartFn$ comes equipped with a natural data service structure (with respect to the Cartesian product of sets as tensor product) such that the inequalities described above hold for all morphisms $f$ in the respective category.
In fact, \eqref{eq:weakhomcomult} is always an equality in $\PartFn$, though it's only an inequality in $\Rel$ and other non-deterministic contexts.
Since we are using $\Poset$-enrichments to describe partial definition of operations, requiring these inequalities in $\Data'(\CC)$ helps us to accurately reflect the desired behavior.

\begin{defn}
A \emph{$\Poset$-enriched attribute structure} on a $\Poset$-enriched SMC $\CC$ is an attribute structure $(\AA, E, V, \gamma)$ such that:
\begin{itemize}
    \item $\AA$, $E$, and $V$ are all $\Poset$-enriched;
    \item $V$ factors through the inclusion $\Data'(\CC) \hookrightarrow \Data(\CC)$.
\end{itemize}
\end{defn}

\subsection{An Example} \label{sec:example}

A key motivating example for our definition of categories with attributes is the modeling of planning in robotics.
This context provides a collection of excellent simple (and not-so-simple) examples of the interaction of physical and informational resources.
We walk through a toy example in this context in the hopes of clarifying our definitions above.
This example will be described in the non-enriched case, though readers can easily extend it to the $\Poset$-enriched case as desired.

It should be noted that the purpose of this example is to provide intuition for the rest of the paper and inspiration for future research into these themes.
As such, this example may need modifications before being used in practical applications or being connected to the functorial semantics below.

\begin{ex}
Suppose that we have a strict SMC $\CC$ with:
\begin{itemize}
    \item generator objects $R$ (robot), $B$ (ball), $R_B$ (robot holding ball), and $L$ (location);
    \item morphisms $\mu, \delta, \epsilon$ such that $(L, \mu, \delta, \epsilon)$ is a data service;
    \item for each generator object $X$ other than $L$, a morphism $\gamma_X$ comprising an action of $(L, \delta, \epsilon)$ on $X$;
    \item and morphisms $\MoveTo: R \otimes L \to R$, $\MoveTo': R_B \otimes L \to R_B$, $\Pick: R \otimes B \to R_B$, and $\Place: R_B \to R \otimes B$.
\end{itemize}

Suppose furthermore that the following tuple $(\AA, E, V, \gamma)$ defines a valid attribute structure on $\CC$.
Let $\AA$ contain objects $X_L$ where $X \in \Ob \CC$ appears as the domain or codomain of one of $\MoveTo$, $\MoveTo'$, $\Pick$, or $\Place$ above.
Furthermore, let $\AA$ be generated by morphisms $f_L: X_L \to Y_L$ where $f$ is one of $\MoveTo, \MoveTo', \Pick,$ or $\Place$, $X$ is the domain of $f$, and $Y$ is the codomain of $f$.
Define $E: \AA \to \CC$ by $E(X_L) = X$ and $E(f_L) = f_L$ for all $X_L$ and $f_L$ in $\AA$.

The functor $V$ is more complicated to define.
Let $V(X_L) = L$ for all $X_L$ such that $X$ is a generator object of $\CC$, and let $V((X \otimes Y)_L) = L \otimes L$ when $X$ and $Y$ are both generator objects of $\CC$ (here the data service structure on $L \otimes L$ is defined in the obvious way).
Define $V$ on the generator morphisms by:
\begin{align}
V(\MoveTo) &= \epsilon \otimes \id_L \\
V(\MoveTo') &= \epsilon \otimes \id_L \\
V(\Pick) &= \mu \\
V(\Place) &= \delta
\end{align}

To define the requisite natural transformation $\gamma$, first note that, given actions $\gamma_X$ and $\gamma_Y$ of the comonoid $L$ on $X$ and $Y$ respectively, we can define an action $\gamma_{X \otimes Y}$ of the comonoid $L \otimes L$ on $X \otimes Y$ by
\begin{equation}\begin{tabular}{cc}
\begin{tabular}{c}
$\gamma_{X \otimes Y} :=$
\end{tabular}
&
\begin{tabular}{c}
    \begin{tikzpicture}
        
    	% Frame
    	\coordinate (NW) at (0,1); \coordinate (NE) at (3,1); 
    	\coordinate (SW) at (0,0); \coordinate (SE) at (3,0);
    	%\draw[black] (NW) rectangle (SE);
    	
    	%Nodes
    	\node[AttrR] (AttrR1) at (0.75, 0.75) {};
    	\node[AttrR] (AttrR2) at (1.35, 0.25) {};
    	
    	% Autocreate port names
    	\pgfmathsetmacro{\diagramIn}{4}
    	\pgfmathsetmacro{\diagramOut}{4}	
    	\pgfmathsetmacro{\nAttrL}{0}
    	\pgfmathsetmacro{\nAttrR}{2}
    	\pgfmathsetmacro{\nMoveTo}{0}
    	\pgfmathsetmacro{\nPick}{0}
    	\pgfmathsetmacro{\nPlace}{0}

    	% Draw edges
    	\path[every node/.style={font=\scriptsize}]
    	(in1) edge[ob] node[at start, yshift=\port] {$X$} (AttrR1in)
        (AttrR1out) edge[ob] node[at end, yshift=-\port] {$X$} (out1)
        (AttrR1) edge[info, out=0] node[at end, yshift=-\port] {$L$} (out3)
    	(in2) edge[ob] node[at start, yshift=\port] {$Y$} (AttrR2in)
    	(AttrR2out) edge[ob] node[at end, yshift=-\port] {$Y$} (out2)
        (AttrR2) edge[info, out=0] node[at end, yshift=-\port] {$L$} (out4)
    	;
    \end{tikzpicture}
\end{tabular}
\end{tabular}\end{equation}
Use this approach to define $\gamma_{R \otimes L}, \gamma_{R_B \otimes L},$ and $\gamma_{R \otimes B}$.
Then, for each $X_L \in \Ob(\AA)$, let $\gamma_{X_L} = \gamma_X$.
This completes the tuple $(\AA, E, V, \gamma)$.
\end{ex}

Viewing $\CC$ as a category with attributes in this way enables us to perform useful reasoning about $\CC$.
For example, given the above attribute structure, we can derive \eqref{eq:moveto} as well as both of the equations shown in Figure  \ref{fig:proofs}.
We can also show that $\Pick \circ \chi_{R \otimes B} = \Pick$ and $\chi_{R \otimes B} \circ \Place = \Place$, where $\chi_{R \otimes B}$ is defined using $\gamma_R$ and $\gamma_B$ as in \eqref{eq:chidef}.

Furthermore, this attribute structure allows us to analogize $\MoveTo$ to the ``put'' of a lens, with $\gamma_R$ serving as the corresponding ``get.''
From this perspective, we can formulate laws on $\MoveTo$ that admit similar interpretations to the Put-Put, Put-Get, and Get-Put lens laws, albeit taking into account the resource- and time-sensitive nature of these operations.

From this viewpoint, the coassociativity property of $\gamma_R$ can be viewed as a novel ``Get-Get'' law, controlling what happens when one retrieves an attribute twice in a row.
Further research is necessary to understand the connections between lenses and resource-sensitive situations like that described in this example.

\section{Semantics}
\label{sec:semantics}

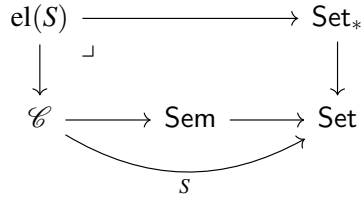
\begin{figure}
    \centering
    \begin{tikzcd}
    \el(S) \ar[rr] \dar \ar[dr, phantom, "\lrcorner", very near start] & & \Set_* \dar \\
    \CC \rar \ar[rr, bend right, "S"'] & \Sem \rar & \Set
    \end{tikzcd}
    \caption{Schematic of the discussion at the beginning of Section \ref{sec:semantics}. Special cases of interest are $\Sem = \BoolAlg^\textrm{op}$ and $\Sem = \Geom$.}
    \label{fig:semantics}
\end{figure}

The semantics of a category with attributes $\CC$ typically manifests itself through a functor from $\CC$ to some ``semantics category'' $\Sem$.
Oftentimes, there also exists a forgetful functor $\Sem \to \Set$, allowing us to define a composite functor $S: \CC \to \Set$.
The category of elements $\el(S)$ of $S$ also reflects information about the semantics of $\CC$; while $\Sem$ can be viewed as a category of state spaces and state space transformations, $\el(S)$ can be viewed as a category of states and state update rules.
This allows us to make mathematical sense of intuitive terms like ``instance of an entity'' that were used above.
See Figure \ref{fig:semantics} for a depiction of this setup.

In this section, we examine the semantics of categories with attributes by way of two key examples (Boolean and geometric semantics) relevant to applications in engineering.
We also illustrate the above construction involving categories of elements in the case of geometric semantics.

\subsection{Boolean Semantics}

The Planning Domain Definition Language (PDDL) is a standard encoding for ``classical'' planning problems involving collections of Boolean variables and discrete actions that operate on them \cite{PDDL}. It extends the earlier STRIPS\footnote{Stanford Research Institute Problem Solver} language with features including typing and object equality. PDDL aligns quite well with the diagrammatic syntax presented earlier, providing both a graphical interface to define and report PDDL problems and solutions, and an approach to string-diagram generation based on existing PDDL solvers.

PDDL can be divided into three fundamental units: \emph{domains}, \emph{problems} and \emph{solutions}. Domains define the relevant (Boolean) features of the situation of interest, as well as the actions that manipulate these features. A problem is posed relative to a domain; it specifies a collection of elements, their initial state, and a desired goal state. A solution gives a sequence of applied actions which transform the initial state into the goal.

Before proceeding we briefly recall a few facts about Boolean algebras \cite{StoneSpaces}, which will provide the substance for our high-level semantics. The free Boolean algebra over a finite set of atomic statements $S$ is given by $\Bool(S):=2^{2^S}$; here we think of a map $s:S\to 2$ as a truth function defining a global state, while a proposition defines a map $q:2^S\to 2$ picking out a subset of global states.

Any element of a Boolean algebra $q\in B$ defines an associated quotient algebra $B\epi B/q$; logically speaking, $B/q$ represents the theory $B$ extended by the axiom $q$. A homomorphism $b:B\to B'$ descends to the quotient $B/q$ if and only if $b(q)=\top$, in which case we write $b\models q$.

A \emph{point} of a Boolean algebra $B$ (also called a valuation) is a homomorphism $b:B\to 2$. For a free, finite algebra $\Bool(S)$ it is easy to show that this is equivalent to a truth function $S\to 2$. Stone duality estabilishes a contravariant relationship between algebra homomorphisms and (continuous) functions between points. Here, where everything is finite, the Stone topologies are discrete and every function between points induces a homomorphism in the opposite direction.

A PDDL domain consists of two main elements, a set of \emph{predicates} $P$ and a set of \emph{actions} $A$. Each predicate has a list of variables which are (optionally) typed from a fixed set or hierarchy $T$, corresponding to an arity function $\arity:P\to \List(T)$. Given a set of typed variables $\tp:X\to T$, the collection of atomic statements $\big\{p(\overline{x})\big\}$ can be identified with the pullback $\Atom(P/X):= P\underset{\mathclap{\List(T)}} {\times} \List(X)$. We write $P(X)$ for the free algebra over $\Atom(P/X)$, yielding a functor $P:\Set/T\to\BoolAlg$.

A PDDL action $a\in A$ consists of a set of typed \emph{parameters} $X_a$ along with pre- and post-conditions $q^a_0,q^a_1\in P(X_a)$. Because $P(X_a)$ is free, this induces an associated homomorphism $a:P(X_a)/q^a_1\to P(X_a)/q^a_0$. This is easiest to see by defining a dual function on points, regarded as truth functions $s:\Atom(P / X_a) \to 2$:
\begin{equation} \label{eq:asdefinition}
	a(s): p(\overline{x})\mapsto \left\{\begin{array}{ll} \neg s(p(\overline{x})) & \textrm{if }p(\overline{x})\wedge q^a_1=\bot\\ s(p(\overline{x}))&\textrm{otherwise}\\\end{array}\right.
\end{equation}
In other words, we flip any bits that are inconsistent with the post-condition, and leave everything else alone.

From the actions $A$, we construct a pair of categories $\BB$ and $\CC$ along with a span of functors $\CC\from \BB \stackrel{I}{\to} \el(\Pt)$, where $\el(\Pt)$ is the category of elements of the functor $\Pt: \BoolAlg^\textrm{op} \to \Set$ sending a Boolean algebra to its set of points. Moreover, the functor $\BB\to\CC$ is a ``partial opfibration,'' in that lifts of arrows are unique if they exist.

Both $\BB$ and $\CC$ are free symmetric monoidal. In $\BB$, the objects are Boolean points, and we include one generator $a_s:s\to a(s)$ for each point $s\models q^a_0$. By contrast, $\CC$ is generated by the objects $X_a$ (for $a \in A$) and contains a generator morphism $a: X_a \to X_a$ for each $a \in A$. The functor $\BB \to \CC$ forgets about the pre- and post-conditions that appear in $\BB$.

There is an obvious functor $I:\BB\to\el(\Pt)$ sending $a_s$ to the point function defined in \eqref{eq:asdefinition}. We note that $I$ does not preserve the monoidal structure because a global state over two collections $X$ and $Y$ will include predicates with variables from both collections. However, $I$ is \emph{oplax} monoidal: a global state over $X$ and $Y$ entails global states on $X$ and $Y$ individually, corresponding to a function $I(X\otimes Y)\to IX\times IY$.

A PDDL problem is defined relative to a domain, and specifies a set of typed variables (called \emph{objects}) $O$. In addition, a problem provides initial and goal states, formulated as propositions $q_0^*,q_1^*\in P(O)$. We can regard these as a pair of (category-theoretic) objects in $\BB$.

Finally, a PDDL solution is given as a sequence of (validly) applied actions which, when applied to any initial state $s_0\models q_0^*$, result in a final state $s_1\models q_1^*$. Here an ``applied action'' for a state $s\in P(O)$ consists of an action $a$ and a function $j:X_a\to O$. This induces a function $j^*:\Pt(P(O))\to\Pt(P(X_a))$, and for an application to be valid we should have $j^*s\models q_0^a$.

Because $P(O)$ is free, we can use the principle of minimal modification to lift the local transformation $a_{j^*s}$ to a global map $a^j_s:\Pt(P(O))\to\Pt(P(O))$ satisfying $j^*(a^j_s(s))\models q^a_1$; i.e., the new global state satisfies the local post-condition for $a$. Given a sequence of applied actions $(a_i,j_i)$ we can iteratively construct a sequence of states $s_i$ and, if each application is valid for the previous state, this will define a string diagram in $\BB$.

We can relate our Boolean semantics and the attribute syntax introduced earlier by associating any free attribute category with a PDDL domain. The types of the domain are the atomic entities and data services of the syntax. Predicates are defined from the attributes. Given two attributes over a shared data service, say with underlying comonoid actions $\phi:M\to M\otimes D$ and $\psi:N\to N\otimes D$, we introduce a binary predicate $D^{\phi,\psi}(M,N)$ indicating whether or not their data values agree. (Based on our semantic intuition from before, $D^{\phi,\psi}(M,N)$ should be true precisely when the $\chi$ morphism associated to these attributes is defined.)

Actions of the PDDL domain correspond to generating morphisms of the syntax, and the equational axioms attached to a morphism define the pre- and post-conditions of the action. 
For example, the equation shown in (\ref{eq:moveto}) corresponds to a post-condition $L^{\phi,\delta}(R,L)$ on the $\MoveTo$ action, equating the location attribute of the robot with the copy ``attribute'' of the target location.

\subsection{Geometric Semantics}

Here we describe a ``geometric semantics category'' $\Geom$ and its connections to categories with attributes and modeling.
The category $\Geom$ provides a lower-level semantic counterpart to the Boolean approach outlined above.
In addition, it yields a ``physical'' approach to understanding the behaviors specified by a category with attributes.
Informally, objects in $\Geom$ are physical objects together with data values, while morphisms in $\Geom$ are (partially defined) paths of the physical objects and update rules for the corresponding values.

Formally, an object $X = (\{ X_i \}, P_X, \theta_X)$ in $\Geom$ consists of a sequence $\{ X_i \}_{i=1}^{k_X}$ of subsets of $\R^3$ where each $X_i$ is called a \emph{simple object}, a topological space $P_X$ called the \emph{parameter space}, and a continuous function $\theta_X = (\theta_{X,i}) : P_X \to \SE(3)^{k_X}$, called the \emph{structure map},\footnote{
Here $\SE(3)$ is the group of rigid motions in $\R^3$.
}
such that for all $p \in P_X$, the sets $\theta_{X,i}(p) \cdot U_i$ are pairwise disjoint.
In general, if $X$ is an object in $\Geom$, we will write $\{ X_i \}_{i=1}^{k_X}$ for the simple objects of $X$, $P_X$ for the parameter space of $X$, and $\theta_X$ for the structure map of $X$.
We use similar notation when the object is instead called $Y$ or anything else.

A good example object to keep in mind is a multi-jointed robot arm.
Such an arm can be modeled as a finite collection of rods, the positions of which vary based on parameters (e.g.\ Euler angles).
The condition imposed on the structure map ensures that, regardless of the parameters, no two rods ever occupy the same point in space at once.

%(here we are using the standard action of $\SE(3)$ on $\PP(\R^3)$, i.e.\ $T \cdot U = \{ T(x) : x \in U \}$).
%\footnote{
%There is no mathematical reason to restrict our attention to $\R^3$ and $\SE(3)$ here; we do so merely for ease of visualization and with an eye towards physical applications.
%}

A morphism $f: X \to Y$ in $\Geom$, where $X$ and $Y$ have the same simple objects (in the same order), consists of a continuous partial function $\Phi_f: P_X \to P_Y$ and a continuous partial function $\phi_f: P_X \times [0, T_f] \to \SE(3)^{k_X}$ (for some $T_f \geq 0$),
such that:
\begin{itemize}
    \item For all $p \in P$ and $t \in [0, T_f]$, $\phi_f(p, t)$ is defined if and only if $\Phi_f(p)$ is defined;
    \item For all $p \in P$ such that $\Phi_f(p)$ is defined:
    \begin{itemize}
        \item $\phi_f(p, 0) = \theta_X(p)$;
        \item $\phi_f(p, T_f) = \theta_Y(\Phi_f(p))$;
        \item For all $t \in [0, T_f]$, the sets $\pi_i(\phi_f(p, t)) \cdot U_i$ are pairwise disjoint.
    \end{itemize}
\end{itemize}
If $X$ and $Y$ are two objects of $\Geom$ such that $X$ and $Y$ do not have the same sequence of simple objects, then there are no morphisms between $X$ and $Y$.

Morphisms $f: X \to Y$ and $g: Y \to Z$, specified using the pattern established above, 
%(with $f$ described by $\Phi_f: P_X \to P_Y$ and $\phi_f: P_X \times [0, T_f] \to \SE(3)^k$ and $g$ described by $\Phi_g: P_X \to P_Y$ and $\phi_g: P_X \times [0, T_g] \to \SE(3)^k$)
can be composed as follows. 
We set $\Phi_{g \circ f} = \Phi_g \circ \Phi_f$ (composition of partial functions) and define $\phi_{g \circ f}$ by:
\[
\phi_{g \circ f}(p, t) = \begin{cases}
\phi_f(p, t) & t \in [0, T_f] \text{ and $\Phi_{g \circ f}(p)$ is defined} \\
\phi_g(\Phi_f(p), t - T_f) & t \in [T_f, T_f + T_g] \text{ and $\Phi_{g \circ f}(p)$ is defined.}
\end{cases}
\]
It is clear that identity morphisms exist and that composition is associative, so $\Geom$ is in fact a category.
In fact, $\Geom$ is a $\Poset$-enriched category, where $f \leq g$ if and only if $T_f = T_g$ and $\Phi_f \leq \Phi_g$ and $\phi_f \leq \phi_g$ as partial functions.

The category $\Geom$ admits a useful notion of ``tensor product of objects,'' defined as follows.
Let $X, Y \in \Ob (\Geom)$; then $X \otimes Y$ is defined to be the object $(\{ X_i \}_{i=1}^{k_X} \sqcup \{ Y_j \}_{j=1}^{k_Y}, P_X \otimes P_Y, (\theta_X \times \theta_Y) |_{P_X \otimes P_Y})$, where $P_X \otimes P_Y$ is the subspace of $P_X \times P_Y$ given by the set:
\[
\{ (p_x, p_y) \in P_X \times P_Y : \forall i, j, \hspace{1em} (\theta_{X,i}(p_x) \cdot X_i) \cap (\theta_{Y,j}(p_y) \cdot Y_j) = \emptyset \}.
\]
In other words, $P_X \otimes P_Y$ is the largest subset $P$ of $P_X \times P_Y$ for which $(\{ U_i \}_{i=1}^{k_X} \cup \{ V_j \}_{j=1}^{k_Y}, P, (\theta_X \times \theta_Y) |_{P})$ is a valid object of $\Geom$.
It is not clear how to extend or adapt this notion to allow for ``tensor products of morphisms,'' which would equip $\Geom$ with a symmetric monoidal structure or something of that sort.\footnote{
The difficulty lies in the interchange law $(f_1 \otimes g_1) \circ (f_2 \otimes g_2) = (f_1 \cdot f_2) \otimes (g_1 \circ g_2)$; when morphisms represent paths (as above), rather than paths up to some suitable notion of homotopy or reparametrization, this law does not hold for ``obvious'' definitions of $\otimes$.
}

We can specify the geometric semantics of a given category with attributes $\CC$ by defining a functor $F: \CC \to \Geom$.
(This $F$ may be $\Poset$-enriched if the attribute structure on $\CC$ is.)
In order to construct such an $F$, we would need models of all the relevant objects and processes in $\CC$.
Hence, instead of outright defining any such $F$ here, we describe properties $F$ should satisfy.

We expect that $F$ sends an entity $A \in \Ob \CC$ to an object $F(A) \in \Geom$ with simple objects consisting of the component parts of some physical model of $A$, parameter space describing the possible states of $A$, and structure map attaching to each state a physical configuration of the component parts.
This approach works even if $A$ is a value; in this case, $F(A)$ has no simple objects ($k_{F(A)} = 0$) and so is described entirely by its parameter space.

The functor $F$ sends a morphism $f: A \to B$ to a physical model of the process specified by $f$.
Specifically, given a state $p \in P_{F(A)}$, the map $\Phi_{F(f)}$ sends $p$ to the state achieved by applying $f$ to a system in the state $p$ (if such a state exists), while $\phi_{F(f)}$ describes the physical motion needed to accomplish this change in state.

When $A$ is a value / data service, the morphisms $\mu_A, \delta_A, \epsilon_A$ making $A$ into a data service would typically be sent to morphisms much like those of the canonical data service structures on objects of $\PartFn$.
For example, $\Phi_{F(\delta_A)}$ might be the partial function $P_A \to P_A \times P_A$ given by $p \mapsto (p, p)$; in this case, we would have $T_{F(\delta_A)} = 0$ and $\phi_{F(\delta_A)}(p, 0) = \Phi_{F(\delta_A)}(p) = (p, p)$.
The image of the natural transformation $\gamma$ from the attribute structure could be defined in a similar way.

We can instantiate the discussion at the beginning of this section in the context of geometric semantics as follows.
We have a canonical functor $\Geom \to \PartFn$ (the category of sets and partial functions) defined by $X \mapsto P_X$ and $f \mapsto \Phi_f$.
Composing this functor with the standard equivalence $\PartFn \simeq \Set_*$ (the category of pointed sets and point-preserving maps) and the forgetful functor $\Set_* \to \Set$ yields a forgetful functor $G: \Geom \to \Set$.

Given $F$ and $G$ as above, the functor $S = F \circ G: \Geom \to \Set$ sends $X \in \Ob \CC$ to the parameter space of some model of $X$ and a morphism to its action on such parameter spaces.
If we consider the state of an object to be entirely determined by a point in its parameter space (and also allow the existence of an ``undefined'' state for each object of $\CC$), then the category of elements $\el(S)$ consists of states of objects in $\CC$ and processes in $\CC$ transforming the domain state into the codomain state.

\section{Conclusion}

In this work, we have presented a categorical interpretation of the notion of an ``attribute'' of an object.
We have also discussed example applications of this interpretation to the categorical modeling of robotics, illustrating how categorical perspectives can be used for the benefit of engineering and other applied fields.

Several interesting questions remain about categorical notions of ``attributes.''
For example, it appears that several of the axioms for $\Poset$-enriched categories with elements break down when considered with respect to subprobabilistic semantics, e.g.\ semantics valued in the category $\mathsf{SRel}$ of \cite{panangaden1998probabilistic}, leading one to wonder how the approach here might be modified for compatibility with such semantics.
In addition, it is sometimes useful to think of time as an attribute of an object, although there are enough differences between time and the other attributes discussed here that it seems that a separate formalism might work better for modeling time.
We hope that future work in this area will be able to resolve some of these questions.

% future work: 
% - refine and explore notions of space- and point-semantics?
% - incorporate time into the syntax side of things?
% - find a way to deal with SRel?

\subsection*{Acknowledgments}

The authors would like to thank Angeline Aguinaldo, Blake Pollard, Fred Proctor, and Eswaran Subrahmanian for helpful comments and discussions.
Thanks are also due to the anonymous reviewers for comments that helped the authors improve and clarify aspects of this paper.

The second named author would also like to thank the organizers of the Applied Category Theory 2020 Adjoint School for providing an excellent environment in which to gain a deeper understanding of category theory. 
In particular, special thanks are due to Paolo Perrone, whose references to categories of elements during meetings of the School helped the second author to realize the role that such constructions played in the semantics of categories with attributes.

\subsection*{Disclaimer}

This paper includes contributions from the U. S. National Institute of Standards and  Technology, and is not subject to copyright in the United States. Commercial products are identified in this article to adequately specify the material. This does not imply recommendation or endorsement by the National Institute of Standards and Technology, nor does it imply the materials identified are necessarily the best available for the purpose.

\bibliographystyle{eptcs}
\bibliography{AttrCats}

\end{document}